%% A manuscript in PLAIN TeX
\def\author{E. Esteves, M. Gagn\'e, and S. Kleiman}
\def\title{Autoduality of the compactified Jacobian}
\def\date{November 9, 1999}
\def\abstract{
 We prove the following {\it autoduality theorem} for an integral
projective curve $C$ in any characteristic.  Given an invertible sheaf
$\cL$ of degree 1, form the corresponding Abel map
  $A_\cL\:C\to\CJ$,
which maps $C$ into its compactified Jacobian, and form its pullback
map $A_\cL^*\:\Pic^0_\CJ\to J$, which carries the connected
component of 0 in the Picard scheme back to the Jacobian.  If $C$ has,
at worst, points of multiplicity 2, then $A_\cL^*$ is an isomorphism,
and forming it commutes with specializing $C$.

Much of our work is valid, more generally, for a family of curves with,
at worst, points of embedding dimension 2.  In this case, we use the
determinant of cohomology to construct a right inverse to
$A_\cL^*$.  Then we prove a scheme-theoretic version of the theorem
of the cube, generalizing Mumford's, and use it to prove that
$A_\cL^*$ is independent of the choice of $\cL$.  Finally, we prove
our autoduality theorem: we use the presentation scheme to achieve an
induction on the difference between the arithmetic and geometric genera;
here, we use a few special properties of points of multiplicity 2.
 }
%%%%%%%%%%%%%%%%%       SOME CHOICES
\def\TheMagstep{\magstep1}      % Normal magnification
                %% Changed to \magstep0 by \DoublepageOutput{TRUE}
\def\PaperSize{letter} % \PaperSize is used to
% \def\PaperSize{AFour}         %  center text on page
        %%%%%%%%%
%%% To get two pages side by side in landscape mode on an ordinary
%%%    sheet via dvips and  PostScript, enable the next command.
% \def\DoublepageOutput{TRUE}
%%%%%%%%%%%%%%%%%       END of CHOICES
%%%%%%%%%%%%%%%%%       Left-Right
 % to put page on right
% \def\FirstPageOnRight{TRUE}

%%%%%%%%%%%%%%%%%       MACRO FILES
%\input acj.mac                          %% Format and style
%%%%%%%%% emd.mac = a Plain TeX macro package by SLK
 \let\@=@  %% For e-mail address

 %% For Repeating Similar Acts Defining Macros
\def\GetNext#1 {\def\NextOne{#1}\if\relax\NextOne\let\next=\relax
        \else\let\next=\DoIt \fi \next}
\def\DoIt{\Act\NextOne\GetNext}
\def\ActOn#1{\expandafter\GetNext #1\relax\ }
\def\defcs#1{\expandafter\xdef\csname#1\endcsname}

 %% PAGE LAYOUT
\parskip=0pt plus 1.8pt \parindent10pt
\hsize 30pc%32pc %29pc %29
\vsize 45pc %48pc %44
\abovedisplayskip 4pt plus3pt minus1pt
\belowdisplayskip=\abovedisplayskip
\abovedisplayshortskip 2.5pt plus2pt minus1pt
\belowdisplayshortskip=\abovedisplayskip

\def\TRUE{TRUE} % For Boolean tests
\ifx\DoublepageOutput\TRUE \def\TheMagstep{\magstep0} \fi
\mag=\TheMagstep

% CENTER TEXT ON PAGE
        % additional vertical adjustment
\newskip\vadjustskip %\vadjustskip=0.5\normalbaselineskip
\def\centertext
 {\hoffset=\pgwidth \advance\hoffset-\hsize
  \advance\hoffset-2truein \divide\hoffset by 2\relax
  \voffset=\pgheight \advance\voffset-\vsize
  \advance\voffset-2truein \divide\voffset by 2\relax
  \advance\voffset\vadjustskip
 }
\newdimen\pgwidth\newdimen\pgheight
\def\letter{letter}\def\AFour{AFour}
\ifx\PaperSize\letter
 \pgwidth=8.5truein \pgheight=11truein
 \message{- Got a paper size of letter.  }\centertext
\fi
\ifx\PaperSize\AFour
 \pgwidth=210truemm \pgheight=297truemm
 \message{- Got a paper size of AFour.  }\centertext
\fi

%% HEADLINE STYLE
\def\today{\ifcase\month\or     % From the TeX book p. 406
 January\or February\or March\or April\or May\or June\or
 July\or August\or September\or October\or November\or December\fi
 \space\number\day, \number\year}
% \nopagenumbers
\footline={\hss\eightpoint\folio\hss}
 \newcount\pagenumber \pagenumber=1
 \def\advancepagenumber{\global\advance\pagenumber by 1}
\def\folio{\number\pagenum} % \pagenum is \let below by an \ifx
\headline={%
  \ifnum\pagenum=0\hfill
  \else
   \ifnum\pagenum=1\firstheadline
   \else
     \ifodd\pagenum\oddheadline
     \else\evenheadline\fi
   \fi
  \fi
}
\expandafter\ifx\csname date\endcsname\relax \let\dato=\today
            \else\let\dato=\date\fi
\let\firstheadline\hfill
\def\oddheadline{\eightpoint
 \rlap{\dato} \hfil \headtitle
 %\hfil\llap{\folio}
 }
\def\evenheadline{\eightpoint %\rlap{\folio} \hfil
 \author\hfil\llap{\dato}}
\def\headtitle{\title}

%% TWO-COLUMN LANDSCAPE FORMAT
% Modified from the TeX book, p. 257.
 \newdimen\fullhsize \newbox\leftcolumn
 \def\fulline{\hbox to \fullhsize}
\def\doublepageoutput
{\let\lr=L
 \output={\if L\lr
           \global\setbox\leftcolumn=\columnbox \global\let\lr=R%
          \else \doubleformat \global\let\lr=L
          \fi
        \ifnum\outputpenalty>-20000 \else\dosupereject\fi
        }%
 \def\doubleformat{\shipout\vbox{%
     \ifx\PaperSize\AFour
           \fulline{\hfil\box\leftcolumn\hfil\columnbox\hfil}%
     \else
           \fulline{\hfil\hfil\box\leftcolumn\hfil\columnbox\hfil\hfil}%
     \fi             }%
     \advancepageno
}
 \def\columnbox{\vbox
   {\if E\topmark\headline={\hfil}\nopagenumbers\fi
    \makeheadline\pagebody\makefootline\advancepagenumber}%
   }%
\fullhsize=\pgheight \hoffset=-1truein
 \voffset=\pgwidth \advance\voffset-\vsize
  \advance\voffset-2truein \divide\voffset by 2
  \advance\voffset\vadjustskip
 
 %%% to put page on right
\ifx\FirstPageOnRight\TRUE % to put page on right
 \null\vfill\nopagenumbers\eject\pagenum=1\relax
\fi
}
\ifx\DoublepageOutput\TRUE \let\pagenum=\pagenumber\doublepageoutput
 \else \let\pagenum=\pageno \fi

%% ADDITIONAL FONTS
 \font\twelvebf=cmbx12          % For title
 \font\smc=cmcsc10              % For authors' names and section headings
%% A LARGE LOWER CASE \OMEGA
%\font\tenbi=cmmi10 scaled \magstep2 % big (math) italic (FROM BEFORE)
 \font\tenbi=cmmi14
 \font\sevenbi=cmmi10 \font\fivebi=cmmi7
 \newfam\bifam  \textfont\bifam=\tenbi
 \scriptfont\bifam=\sevenbi \scriptscriptfont\bifam=\fivebi
 \mathchardef\variablemega="7121 
 \mathchardef\variablenu="7117 
%% EIGHT POINT TYPE for footnotes and references
\catcode`\@=11          % make @ a letter temporarily
\def\eightpoint{\eightpointfonts
 \setbox\strutbox\hbox{\vrule height7\p@ depth2\p@ width\z@}%
 \eightpointparameters\eightpointfamilies
 \normalbaselines\rm
 }
\def\eightpointparameters{%
 \normalbaselineskip9\p@
 \abovedisplayskip9\p@ plus2.4\p@ minus6.2\p@
 \belowdisplayskip9\p@ plus2.4\p@ minus6.2\p@
 \abovedisplayshortskip\z@ plus2.4\p@
 \belowdisplayshortskip5.6\p@ plus2.4\p@ minus3.2\p@
 }
\newfam\smcfam
\def\eightpointfonts{%
 \font\eightrm=cmr8 \font\sixrm=cmr6
 \font\eightbf=cmbx8 \font\sixbf=cmbx6
 \font\eightit=cmti8
 \font\eightsmc=cmcsc8
 \font\eighti=cmmi8 \font\sixi=cmmi6
 \font\eightsy=cmsy8 \font\sixsy=cmsy6
 \font\eightsl=cmsl8 \font\eighttt=cmtt8}
\def\eightpointfamilies{%
 \textfont\z@\eightrm \scriptfont\z@\sixrm  \scriptscriptfont\z@\fiverm
 \textfont\@ne\eighti \scriptfont\@ne\sixi  \scriptscriptfont\@ne\fivei
 \textfont\tw@\eightsy \scriptfont\tw@\sixsy \scriptscriptfont\tw@\fivesy
 \textfont\thr@@\tenex \scriptfont\thr@@\tenex\scriptscriptfont\thr@@\tenex
 \textfont\itfam\eightit        \def\it{\fam\itfam\eightit}%
 \textfont\slfam\eightsl        \def\sl{\fam\slfam\eightsl}%
 \textfont\ttfam\eighttt        \def\tt{\fam\ttfam\eighttt}%
 \textfont\smcfam\eightsmc      \def\smc{\fam\smcfam\eightsmc}%
 \textfont\bffam\eightbf \scriptfont\bffam\sixbf
   \scriptscriptfont\bffam\fivebf       \def\bf{\fam\bffam\eightbf}%
 \def\rm{\fam0\eightrm}%
% \tt \ttglue=0.5em plus0.25em minus0.15em
 }
%% Modification of the PLAIN footnote macro for 8pt
\def\vfootnote#1{\insert\footins\bgroup
 \eightpoint\catcode`\^^M=5\leftskip=0pt\rightskip=\leftskip%% only change
 \interlinepenalty\interfootnotelinepenalty
  \splittopskip\ht\strutbox % top baseline for broken footnotes
  \splitmaxdepth\dp\strutbox \floatingpenalty\@MM
  \leftskip\z@skip \rightskip\z@skip \spaceskip\z@skip \xspaceskip\z@skip
  \textindent{#1}\footstrut\futurelet\next\fo@t}

%%  ``Ties'' with a \thinspace for page numbers
\def\p.{p.\penalty\@M \thinspace}
\def\pp.{pp.\penalty\@M \thinspace}
%% SECTIONING
\newcount\sctno \newskip\sctnskip \sctnskip=0pt plus\baselineskip
\def\sctn#1\par
  {\removelastskip\vskip\sctnskip %\penalty-50%-250
  \vskip-\sctnskip \bigskip\medskip%\smallskip%\bigskip
  \centerline{#1}\nobreak\medskip
}

\def\sct#1 {\sctno=#1\relax\sctn#1. }

%%  STYLE MACROS  %%
%% Redefine \item to give greater indentation than AMSTeX
%   and the roman font within parentheses.
\def\item#1 {\par\indent\indent\indent%\indent
 \hangindent3\parindent
 \llap{\rm (#1)\enspace}\ignorespaces}
%% Define a similar macro without the hanging indentation for assertions
%% and that starts each part with an ordinary \parindent
 \def\inpart#1 {{\rm (#1)\enspace}\ignorespaces}
 \def\part {\par\inpart}

%% ARTICLES
\def\Cs#1){\(\number\sctno.#1)}
\def\part#1 {\par\(#1)\enspace\ignorespaces}

\def\dsc#1 #2.{\medbreak\Cs#1) {\it #2.} \ignorespaces}
%% For setting results
\def\proclaim#1 #2 {\medbreak
  {\bf#1 (\number\sctno.#2).  }\ignorespaces\bgroup\it}
\def\endproclaim{\par\egroup\medskip}
\def\pf{\endproclaim{\bf Proof.} \ignorespaces}
\def\lem{\proclaim Lemma } \def\prp{\proclaim Proposition }
\def\cor{\proclaim Corollary }  
\def\dfn#1 {\medbreak {\bf Definition (\number\sctno.#1).  }\ignorespaces}
\def\rmk#1 {\medbreak {\bf Remark (\number\sctno.#1).  }\ignorespaces}
%% REFERENCING
        % to introduce the keys in order
 \newcount\refno \refno=0        \def\NoKey{*!*}
 \def\MakeKey{\advance\refno by 1 \expandafter\xdef
  \csname\TheKey\endcsname{{\number\refno}}\NextKey}
 \def\NextKey#1 {\def\TheKey{#1}\ifx\TheKey\NoKey\let\next\relax
  \else\let\next\MakeKey \fi \next}
 \def\RefKeys #1\endRefKeys{\expandafter\NextKey #1 *!* }
 \def\SetRef#1 #2,{\hang\llap
  {[\csname#1\endcsname]\enspace}{\smc #2},}
 \newbox\keybox \setbox\keybox=\hbox{[25]\enspace}
 \newdimen\keyindent \keyindent=\wd\keybox
\def\references{\kern-\medskipamount
  \sctn References\par
  \vskip-\medskipamount
  \bgroup   \frenchspacing   \eightpoint
   \parindent=\keyindent  \parskip=\smallskipamount
   \everypar={\SetRef}\par}
\def\endreferences{\egroup}

%% SERIALS
 \def\serial#1#2{\expandafter\def\csname#1\endcsname ##1 ##2 ##3
        {\unskip\ {\it #2\/} {\bf##1} (##2), ##3.}} % \serial{}{}

%% modified \cite code from AMSTeX
\def\UThin{\penalty\@M \thinspace\ignorespaces}
        % unbreakable \thinspace for use after periods
\let\ts=\UThin
\def\(#1){{\let~=\UThin\rm(#1)}}
\def\relaxnext@{\let\next\relax}
\def\cite#1{\relaxnext@
 \def\nextiii@##1,##2\end@{\unskip\space{\rm[\SetKey{##1},\let~=\UThin##2]}}%
 \in@,{#1}\ifin@\def\next{\nextiii@#1\end@}\else
 \def\next{{\rm[\SetKey{#1}]}}\fi\next}
\newif\ifin@
\def\in@#1#2{\def\in@@##1#1##2##3\in@@
 {\ifx\in@##2\in@false\else\in@true\fi}%
 \in@@#2#1\in@\in@@}
\def\SetKey#1{{\bf\csname#1\endcsname}}

\catcode`\@=12  %at signs are no longer letters, but other

 %%  MATH MACROS
\let\:=\colon \let\ox=\otimes \let\x=\times
 
\def\dr#1{#1^\dagger}
\let\into=\hookrightarrow \def\onto{\to\mathrel{\mkern-15mu}\to}
\let\To=\longrightarrow \def\TO#1{\buildrel#1\over\To}
\def\smashedlongrightarrow{\setbox0=\hbox{$\longrightarrow$}\ht0=1pt\box0}
\def\risom{\buildrel\sim\over{\smashedlongrightarrow}}
 \def\lgto{-\mathrel{\mkern-10mu}\to}
 \def\smashedlgto{\setbox0=\hbox{$\scriptstyle\lgto$}\ht0=1.85pt
        \lower1.25pt\box0}
\def\tto{\buildrel\lgto\over{\smashedlgto}}
\let\vf=\varphi  \let\?=\overline

\def\CJ{{\bar J}} \def\smC{C^{\sm}} \def\dM{{\cal M}^\diamond}

 \def\IP{{\bf P}}  
\def\Act#1{\defcs{c#1}{{\cal#1}}}
 \ActOn{A B C D E F G H I J K L M N O P Q R T }
\def\Act#1{\defcs{#1}{\mathop{\rm#1}\nolimits}}
 \ActOn{cod Div Ext Hilb length Pic Quot sm Spec Supp }
\def\Act#1{\defcs{c#1}{\mathop{\it#1}\nolimits}}
 \ActOn{Cok Ext Hom Ker Sym }
\def\Act#1{\defcs{#1}{\hbox{\rm\ #1 }}}
 \ActOn{and by for where with on }

%%% End of the macro file
%\input ../macro/btcd.sty                %% Commutative diagrams
\catcode`\@=11

 \def\activeat#1{\csname @#1\endcsname}
 \def\def@#1{\expandafter\def\csname @#1\endcsname}
 {\catcode`\@=\active \gdef@{\activeat}}

\let\ssize\scriptstyle
\newdimen\ex@   \ex@.2326ex

 \def\requalfill{\cleaders\hbox{$\mkern-2mu\mathord=\mkern-2mu$}\hfill
  \mkern-6mu\mathord=$}
 \def\eqfill{$\m@th\mathord=\mkern-6mu\requalfill}
 \def\deffill{\hbox{$:=$}$\m@th\mkern-6mu\requalfill}
 \def\fiberbox{\hbox{$\vcenter{\hrule\hbox{\vrule\kern1ex
     \vbox{\kern1.2ex}\vrule}\hrule}$}}
 \def\Fiberbox{\rlap{\kern-0.75pt\raise0.75pt\fiberbox}%
        \kern.75pt{\lower0.75pt\fiberbox}}

 \font\arrfont=line10
 \def\Swarrow{\vcenter{\hbox{$\swarrow$\kern-.26ex
    \raise1.5ex\hbox{\arrfont\char'000}}}}

 \newdimen\arrwd
 \newdimen\minCDarrwd \minCDarrwd=2.5pc
        
 \def\findarrwd#1#2#3{\arrwd=#3%
  \setbox\z@\hbox{$\ssize\;{#1}\;\;$}%
 \setbox\@ne\hbox{$\ssize\;{#2}\;\;$}%
  \ifdim\wd\z@>\arrwd \arrwd=\wd\z@\fi
  \ifdim\wd\@ne>\arrwd \arrwd=\wd\@ne\fi}
 \newdimen\arrowsp\arrowsp=0.375em
 \def\findCDarrwd#1#2{\findarrwd{#1}{#2}{\minCDarrwd}
    \advance\arrwd by 2\arrowsp}
 \newdimen\minarrwd 
 \setbox\z@\hbox{$\longrightarrow$} \minarrwd=\wd\z@

 \def\harrow#1#2#3#4{{\minarrwd=#1\minarrwd%
   \findarrwd{#2}{#3}{\minarrwd}\kern\arrowsp
    \mathrel{\mathop{\hbox to\arrwd{#4}}\limits^{#2}_{#3}}\kern\arrowsp}}
 \def@]#1>#2>#3>{\harrow{#1}{#2}{#3}\rightarrowfill}
 \def@>#1>#2>{\harrow1{#1}{#2}\rightarrowfill}
 \def@<#1<#2<{\harrow1{#1}{#2}\leftarrowfill}
 \def@={\harrow1{}{}\eqfill}
 \def@:#1={\harrow1{}{}\deffill}
 \def@ N#1N#2N{\vCDarrow{#1}{#2}\UpDownarrow}
 \def\UpDownarrow{\uparrow\,\Big\downarrow}

\def@'#1'#2'{\harrow1{#1}{#2}\tarrowfill}
 \def\lgTo{\dimen0=\arrwd \advance\dimen0-2\arrowsp
        \hbox to\dimen0{\rightarrowfill}}
 \def\smashedlgTo{\setbox0=\hbox{$\scriptstyle\lgTo$}\ht0=1.85pt
        \lower1.25pt\box0}
 \def\tto{\buildrel\lgTo\over{\smashedlgTo}}
 \def\tarrowfill{\hfil$\tto$\hfil}  % provisional

 \def@={\ifodd\row\harrow1{}{}\eqfill
   \else\vCDarrow{}{}\Vert\fi}
 \def@.{\ifodd\row\relax\harrow1{}{}\hfill
   \else\vCDarrow{}{}.\fi}
 \def@|{\vCDarrow{}{}\Vert}
 \def@ V#1V#2V{\vCDarrow{#1}{#2}\downarrow}
\def@ A#1A#2A{\vCDarrow{#1}{#2}\uparrow}
 \def@(#1){\arrwd=\csname col\the\col\endcsname\relax
   \hbox to 0pt{\hbox to \arrwd{\hss$\vcenter{\hbox{$#1$}}$\hss}\hss}}

 \def\squash#1{\setbox\z@=\hbox{$#1$}\finsm@@sh}
\def\finsm@@sh{\ifnum\row>1\ht\z@\z@\fi \dp\z@\z@ \box\z@}

 \newcount\row \newcount\col \newcount\numcol \newcount\arrspan
 \newdimen\vrtxhalfwd  \newbox\tempbox

 \def\innernewdimen{\alloc@1\dimen\dimendef\insc@unt}
 \def\measureinit{\col=1\vrtxhalfwd=0pt\arrspan=1\arrwd=0pt
   \setbox\tempbox=\hbox\bgroup$}
 \def\setinit{\col=1\hbox\bgroup$\ifodd\row
   \kern\csname col1\endcsname
   \kern-\csname row\the\row col1\endcsname\fi}
 \def\findvrtxhalfsum{$\egroup
  \expandafter\innernewdimen\csname row\the\row col\the\col\endcsname
  \global\csname row\the\row col\the\col\endcsname=\vrtxhalfwd
  \vrtxhalfwd=0.5\wd\tempbox
%  \global\expandafter\advance\csname row\the\row col\the\col\endcsname
  \global\advance\csname row\the\row col\the\col\endcsname by \vrtxhalfwd
  \advance\arrwd by \csname row\the\row col\the\col\endcsname
  \divide\arrwd by \arrspan
  \loop\ifnum\col>\numcol \numcol=\col%
 \expandafter\innernewdimen \csname col\the\col\endcsname
     \global\csname col\the\col\endcsname=\arrwd
   \else \ifdim\arrwd >\csname col\the\col\endcsname
      \global\csname col\the\col\endcsname=\arrwd\fi\fi
   \advance\arrspan by -1 %
   \ifnum\arrspan>0 \repeat}
 \def\setCDarrow#1#2#3#4{\advance\col by 1 \arrspan=#1
    \arrwd= -\csname row\the\row col\the\col\endcsname\relax
    \loop\advance\arrwd by \csname col\the\col\endcsname
     \ifnum\arrspan>1 \advance\col by 1 \advance\arrspan by -1%
     \repeat
    \squash{\mathop{
     \hbox to\arrwd{\kern\arrowsp#4\kern\arrowsp}}\limits^{#2}_{#3}}}
 \def\measureCDarrow#1#2#3#4{\findvrtxhalfsum\advance\col by 1%
   \arrspan=#1\findCDarrwd{#2}{#3}%
    \setbox\tempbox=\hbox\bgroup$}

%$

\def\vCDarrow#1#2#3{\kern\csname col\the\col\endcsname
    \hbox to 0pt{\hss$\vcenter{\llap{$\ssize#1$}}%
     \Big#3\vcenter{\rlap{$\ssize#2$}}$\hss}\advance\col by 1}

 \def\setCD{\def\harrow{\setCDarrow}%
  \def\\{$\egroup\advance\row by 1\setinit}
  \m@th\lineskip3\ex@\lineskiplimit3\ex@ \row=1\setinit}
 \def\endsetCD{$\strut\egroup}
 \def\measure{\bgroup
  \def\harrow{\measureCDarrow}%
  \def\\##1\\{\findvrtxhalfsum\advance\row by 2 \measureinit}%
  \row=1\numcol=0\measureinit}
 \def\endmeasure{\findvrtxhalfsum\egroup}

\newbox\CDbox \newdimen\sdim

 \newcount\savedcount
 \def\CD#1\endCD{\savedcount=\count11%
   \measure#1\endmeasure
   \vcenter{\setCD#1\endsetCD}%
   \global\count11=\savedcount}

 \catcode`\@=\active
%%%%%%%%%%%%%%%%%       REFERENCE KEYS
 \RefKeys
 AIK76 AK76 AK79 AK79II AK80 AK90 E95 E99 EGK FGA EGA K87 KM76 L59 M60 M65
M70 R67
 \endRefKeys
%%%%%%%%%%%%%%%%%       TOPMATTER
{\leftskip=0pt plus1fill \rightskip=\leftskip
 \obeylines
 \leavevmode \bigskip%\bigskip
 {\twelvebf \title
 } \medskip%\bigskip
%% Subject classification and acknowledgements
 \footnote{}{\noindent %
 MSC-class: 14H40 (Primary) 14K30, 14H20  (Secondary).}
%% Authors' names, addresses and support
 Eduardo Esteves\footnote{$^{1}$}{%
    Supported in part by PRONEX and CNPq Proc. 300004/95-8 (NV).}
 {\eightpoint\it\medskip%\bigskip
 Instituto de Matem\'atica Pura e Aplicada
 Estrada D. Castorina {\sl110, 22460--320} Rio de Janeiro RJ, BRAZIL
 \rm E-mail: \tt Esteves\@impa.br \medskip%\bigskip
 }
 Mathieu Gagn\'e
 {\eightpoint\it\medskip%\bigskip
 EMC Corporation
 171 South Street, Hopkinton, MA {\sl01748}, USA
  \rm E-mail: \tt MGagne\@emc.com \medskip%\bigskip
 } and \medskip%\bigskip
 Steven Kleiman\footnote{$^{2}$}{%
    Supported in part by NSF grant 9400918-DMS.}
 {\eightpoint\it\medskip%\bigskip
 Department of Mathematics, Room {\sl 2-278} MIT,
 {\sl77} Mass Ave, Cambridge, MA {\sl02139-4307}, USA
 \rm E-mail: \tt Kleiman\@math.mit.edu \medskip%\bigskip
 \rm \dato \bigskip%\bigskip
 }
}
 %% Abstract
{%\parindent=1.5\parindent \narrower
 \advance\leftskip by 1.5\parindent \advance\rightskip by 1.5\parindent
 \eightpoint \noindent
 {\smc Abstract.}\enspace \ignorespaces \abstract \par}
%%\end of Topmatter

 %% Body
%\input ac1.tex
\sct1 Introduction

Let $C$ be an integral projective curve, defined over an algebraically
closed field of any characteristic, and $\cL$ an invertible sheaf of
degree 1.  Form the (generalized) Jacobian, the connected component of
the identity of the Picard scheme, $J:=\Pic^0_C$.  If $C$ is smooth,
then $J$ is an Abelian variety, and the Abel map
 $A_\cL\:C\to J$
 is defined by $P\mapsto\cL(-P)$.  Also, the corresponding pullback is an
isomorphism, $A_\cL^*\:\Pic^0_J \risom J$, which is independent of
the choice of $\cL$; thus $J$ is ``autodual,'' or canonically isomorphic
to its own dual Abelian variety $\Pic^0_J$.  (See Theorem 3 on \p.156 in
\cite{L59} or Proposition 6.9 on \p.118 in \cite{M65}.) Our main
  result is the autoduality theorem of (2.1); it
 asserts that, more generally, if $C$ has, at worst, double points
(arbitrary points of multiplicity 2), then a similar pullback is an
isomorphism, and forming it commutes with specializing $C$.

Suppose first that $C$ has arbitrary singularities.  Recall (see
\cite{AK76}, \cite{AK79II}, \cite{AK80}) that $J$ has a natural
compactification $\CJ$, the (fine) moduli space of torsion-free sheaves
of rank 1 and degree 0.  Also, the Abel map
 $A_\cL\:C\to \CJ$ is
 defined by
 $P\mapsto\cI_P\ox\cL$
 where $\cI_P$ is the ideal of $P$; it is a closed embedding if $C$ is
not of genus 0.  Furthermore, the Picard scheme $\Pic_{\CJ}$ exists and
is a union of quasi-projective, open and closed
 subschemes\UThin---\thinspace including $\Pic^0_{\CJ}$ and
$\Pic^\tau_{\CJ}$, which are the connected component of $0$ and the
subscheme of points with multiples in $\Pic^0_{\CJ}$.

Suppose now that all the singularities of $C$ are {\it surficial}, that
is, of embedding dimension 2.  Recall (see \cite{AIK76}) that $\CJ$ is
rather nice; it is a local complete intersection, and is integral and
projective.  So forming $\Pic_{\CJ}$ commutes with specializing $C$, but
conceivably forming $\Pic^0_{\CJ}$ does not.  Nevertheless, we prove two
general results, Propositions (2.2) and (3.7).  The former asserts that
$A_\cL^*\:\Pic^0_\CJ\to J$ has a natural right inverse $\beta$,
which is independent of the choice of $\cL$.  The latter is much deeper,
and asserts that $A_\cL^*$ is itself independent of the choice of
$\cL$.

Suppose finally that all the singularities of $C$ are double points.
Then $A_\cL^*$ is an isomorphism and
 $\Pic^0_{\CJ}=\Pic^\tau_{\CJ}$ by
 our autoduality theorem.
  Now, double points are surficial; hence $\CJ$ is integral.  Moreover,
there exists a scheme parameterizing the torsion-free rank-1 sheaves on
$\CJ$; it is a fine moduli space, and its connected components are
projective.  Let $\?U$ denote the closure of $\Pic^0_\CJ$.  Then the
isomorphism $A_\cL^*\:\Pic^0_\CJ\risom J$ extends to a map
 $\eta\:\?U\to\CJ$
 by Corollary~(4.4); this is our
deepest result, and rests on everything preceding it.  Is $\eta$ an
isomorphism?  Perhaps yes, perhaps no; our work does not appear to
suggest which.

All four of our results are compatible with specializing $C$.  More
precisely, we prove relative versions of them for flat, projective
families of geometrically integral curves over an arbitrary locally
Noetherian base scheme.  Some fibers may be smooth, others not.  A node
may degenerate into a cusp; two nodes may coalesce into a tac.

What happens when all the singularities of $C$ are {\it surficial}?  Is
$A_\cL^*$ an isomorphism then too?  The evidence is mixed.  On the
one hand, Propositions (2.2) and (3.7) suggest so, as they assert that
$A_\cL^*$ has a right inverse $\beta$, and both maps are
independent of the choice of $\cL$.  On the other hand, our proof of
 the autoduality theorem
 suggests not; it doesn't simply fail when $C$ has
singularities of higher multiplicity, rather it suggests that then there
may be a counterexample.

Indeed, to prove
 the autoduality theorem,
 we proceed basically as follows.  We
form $J^1:=\Pic^1_C$, the component of the Picard scheme that
parameterizes the invertible sheaves of degree 1.  Then we put together
the Abel maps $A_\cL$, as $\cL$ varies, to form the Abel map of
bidegree (1,1):
        $$A\:C\x J^1\to \CJ.$$
 This map is studied in the authors' paper \cite{EGK}
 (where, however, the two factors are taken in the opposite order; that
is, $A$ maps $J^1\x C$ into $\CJ$).  In \cite{EGK},
 the following facts are proved.
Suppose that $C$ is Gorenstein.  Then
$A$ is smooth; so its image $V$ is open.  Furthermore, if
$g$ denotes the arithmetic genus, then the complement $\CJ-V$ is of
dimension at most $g-2$ if and only if all the singularities of $C$
are double points.

Hence, if $C$ has higher surficial singularities, then $\CJ$ is
irreducible, and $\CJ-V$ contains a set of codimension 1.  This set
could support a Cartier divisor $D$.  If $D$ exists, then
$A_\cL^*\cO(D)$ is trivial for any $\cL$.  Furthermore, $D$ could
vary in an algebraic family with support on $\CJ-V$ and with two
linearly inequivalent members.  If so,
 then
 their difference would correspond
to a point of $\Pic^0_\CJ$, other than $0$.  Thus $A_\cL^*$ would
not be injective, and we'd have a counterexample.

On the other hand, suppose that all the singularities of $C$ are double
points.  Then $\CJ-V$ is small.  Moreover, $\CJ$ is a local complete
intersection.  So we may
 (and will)
 prove that $A_\cL^*$ is injective
basically as follows.  Let $\cN$ be an invertible sheaf on $\CJ$.  Then
$\cN$ is trivial if its restriction $\cN|V$ is trivial.  In turn, to
show that $\cN|V$ is trivial, we may use descent theory since
$A$ is smooth, so flat.

Suppose that $\cN$ corresponds to a point of $\Pic^0_\CJ$.  Then there
are invertible sheaves $\cN_1$ on $C$ and $\cN_2$ on $J^1$ such that
        $$A^*\cN=p_1^*\cN_1\ox p_2^*\cN_2,$$
 where the $p_i$ are the projections; the existence of the $\cN_i$
results from our general theory of the theorem of the cube, especially
Part (2) of Lemma~(3.6).  Consequently, $A_\cL^*\cN$ is equal to
$\cN_1$, so is independent of the choice of $\cL$, as was asserted
above.

Suppose also that $A_\cL^*\cN$ is trivial for some $\cL$.  We have
to prove that $\cN$ is trivial too.  To begin, note that $\cN_1$ is
trivial, so $A^*\cN$ is equal to $p_2^*\cN_2$.  We proceed by
induction on the arithmetic genus $g$.  Choose a double point $Q$ on
$C$, and blow $Q$ up, getting $\vf\:\dr C\to C$.  Then $\dr C$ too has
only double points by Lemma~(6.4) of \cite{EGK}, and its
arithmetic genus is $g-1$ by Proposition (6.1) of \cite{EGK}.

To relate the compactified Jacobians $\CJ_C$ and $\CJ_{\dr C}$ of $C$
and $\dr C$, we use the presentation scheme $P$ and the maps
$\kappa\:P\to\CJ_C$ and $\pi\:P\to\CJ_{\dr C}$; they are studied in
\cite{EGK}.  By Theorem~(6.3) of \cite{EGK}, because $Q$ is a double
point, $\pi$ is a locally trivial $\IP^1$-bundle.  On each $\IP^1$, the
restriction of $\kappa^*\cN$ is trivial because $\cN$ corresponds to a
point of $\Pic^0_{\CJ_C}$.  Hence, $\kappa^*\cN$ is the pullback of a
sheaf $\dr\cN$ on $\CJ_{\dr C}$.

Because $C$ and $\dr C$ are Gorenstein, there are two natural
commutative diagrams
 $$\CD
        \dr C\x J_{C}^1       @>\Lambda>>             P       \\
           @V 1\x\vf^*VV                         @V\pi VV        \\
         \dr C\x J_{\dr C}^1  @>A_{\dr C}>> \CJ_{\dr C}
  \endCD\qquad\CD
        \dr C\x J_{C}^1  @>\Lambda>>          P       \\
           @V\vf\x1VV                    @V\kappa VV     \\
          C\x J_{C}^1    @>A_C>> \CJ_{C}
  \endCD$$
  by Corollary (5.5) in \cite{EGK}; here $A_C$ and $A_{\dr C}$ are the
Abel maps of $C$ and $\dr C$.  The diagrams imply that, if
$\dr\cL:=\vf^*\cL$, then $A_{\dr\cL}^*\dr\cN$ is equal to
$\vf^*A_\cL^*\cN$, which is trivial by hypothesis.

Since autoduality holds for $\dr C$ by induction, $\dr\cN$ is trivial.
Since the pullback of $\dr\cN$ to $P$
is equal to $\kappa^*\cN$, the latter is
trivial.  So thanks to the commutativity of the second diagram above,
$(\vf\x1)^*A_C^*\cN$ is trivial.  Since $A_C^*\cN$ is
equal to $p_2^*\cN_2$, it follows that $\cN_2$  is trivial.

Since $C$ and $\dr C$ are Gorenstein, by Corollary (5.5) in \cite{EGK},
the second diagram above is Cartesian.  Consider the descent data on
$(\vf\x1)^*A_C^*\cN$ with respect to $\Lambda$; since $\kappa^*\cN$ is
trivial, this data is trivial.  Hence, so is that on $A_C^*\cN$ with
respect to $A_C$ because $\kappa$ is birational.  Therefore $\cN$ is
trivial.  Thus $A_\cL^*$ is injective.

We construct the right inverse $\beta\:J\to\Pic^0_\CJ$ to $A_\cL^*$ by
using the determinant of cohomology $\cD$ along the projection
$q_2\:C\x\CJ\to\CJ$.  We proceed as follows.  Fix a universal sheaf
$\cI$ on $C\x\CJ$.  Then, given any invertible sheaf $\cM$ on $C$ of
degree~0, set
        $$\beta(\cM):=(\cD(\cI\ox q_1^*\cM))^{-1}\ox\cD(\cI),$$
 where $q_1\:C\x\CJ\to C$ is the projection.

This construction was suggested by Breen [pvt.\ comm., 1985].  It is a
modern formulation of an older construction using the theta divisor.
Namely,
        $$\beta(\cM)=\Theta_\cL-\tau^*_\cM\Theta_\cL$$
 where $\tau_\cM\:\CJ\to\CJ$ is the translation, given by tensoring
with $\cM$, and where $\Theta_\cL$ is the divisor obtained by pulling
back the canonical theta divisor along the isomorphism
$\CJ\risom\CJ^{g-1}$ given by tensoring with $\cL^{g-1}$.  We consider
the equivalence of the two formulations in more detail in Remark (2.4).

Since $A_\cL^*\beta=1$, the map $\beta$ is a closed embedding.
Since $A_\cL^*$ is injective, we could conclude that it is an
isomorphism if we knew, a priori, that $\Pic^0_\CJ$ is reduced.  We
don't.  So we must prove that $A_\cL^*$ is a monomorphism, that is,
injective on $T$-points; we take care to do so in (4.1).

In short, in Section 2, we formulate the autoduality theorem, our main
result: if the curves in a family have double points at worst, then the
Abel map $A_\cL^*$ is an isomorphism.  Then we treat $\beta$, which
is the canonical right inverse to $A_\cL^*$.  In Section 3, we
generalize Mumford's scheme-theoretic theorem of the cube, and conclude
that $A_\cL^*$ is independent of the choice of $\cL$.  Finally, in
Section 4, we prove our autoduality theorem, and then extend
$A_\cL^*$ to a map from the natural compactification of
$\Pic^0_\CJ$ onto $\CJ$.

\sct2 Autoduality

 \dsc1 Statement.  Consider a {\it flat projective family of integral
curves} $p\:C\to S$; that is, $S$ is a locally Noetherian scheme, and
$p$ is a flat and projective map with geometrically integral fibers of
dimension 1.  Recall (see \cite{AK76}, \cite{AK79II}, \cite{AK80}) that,
given an integer $n$, there exists a projective $S$-scheme $\CJ^n_{C/S}$ that
parameterizes the torsion-free rank-1 sheaves of degree $n$ on the
fibers of $C/S$.  Furthermore, there exists an open subscheme $J^n_{C/S}$
parameterizing those sheaves that are invertible.  Also, forming
$\CJ^n_{C/S}$ and $J^n_{C/S}$ commutes with changing the base $S$.
As is customary, call $J^n_{C/S}$ the (relative generalized) {\it Jacobian}
of $C/S$, and $\CJ^n_{C/S}$ the {\it compactified Jacobian}.  We will
often abbreviate $J^n_{C/S}$ by $J^n$ and $\CJ^n_{C/S}$ by $\CJ^n$. Set
        $$J_{C/S}:=J^0_{C/S}\and\CJ_{C/S}:=\CJ^0_{C/S}.$$
 We will also abbreviate $J_{C/S}$ by $J$ and $\CJ_{C/S}$ by $\CJ$.

More precisely, a (relative) {\it torsion-free rank-1 sheaf $\cI$} on
$C/S$ is an $S$-flat coherent $\cO_C$-module $\cI$ such that, for each
point $s$ of $S$, the fiber $\cI(s)$ is a torsion-free rank-1 sheaf on
the fiber $C(s)$.  Moreover, $\cI$ is {\it of degree} $n$ if $\cI(s)$
satisfies the relation,
        $$\chi(\cI(s))-\chi(\cO_{C(s)})=n.$$

 Given a locally Noetherian $S$-scheme $T$, a torsion-free rank-1 sheaf
of degree $n$ on $C\x T/T$ defines an $S$-map $T\to \CJ^n$.  Conversely,
every such $S$-map arises from such a sheaf, which is determined up to
tensor product with the pullback of an invertible sheaf on $T$, at least
if the smooth locus of $C/S$ admits a section.  If so, then in
particular the identity map $1_{\CJ^n}$ arises from such a sheaf on $C\x
\CJ^n/\CJ^n$; the latter sheaf is known as a {\it universal\/} (or
Poincar\'e) sheaf, as any $T\to\CJ^n$ arises from the sheaf on $C\x
T/T$ obtained by pulling back a universal sheaf.

In general, an $S$-map $T\to\CJ^n$ arises rather from a pair $(T'/T,\,
\cI')$ where $T'/T$ is an \'etale covering (that is, the map $T'\to T$
is \'etale, surjective, and of finite type) and where $\cI'$ is a
torsion-free rank-1 sheaf of degree $n$ on $C\x T'/T'$.  Such a pair
defines such an $S$-map if and only if there is an \'etale covering
$T''\big/\, T'\x_TT'$ such that the two pullbacks of $\cI'$ to $C\x T''$
are equal.  A second such pair $(T'_1/T,\,\cI'_1)$ defines the same
$S$-map if and only if there is an \'etale covering $T''\big/\,
T'\x_TT'_1$ such that the pullbacks of $\cI'$ and $\cI'_1$ to $C\x T''$
are equal.  In sum, $\CJ^n$ represents the \'etale sheaf associated to
the functor of torsion-free rank-1 sheaves.

Given an invertible sheaf $\cL$ of degree 1 on $C/S$, define the {\it
Abel map,}
        $$A_\cL\:C\to \CJ,$$
 as follows.  Let $\cI_{\Delta}$ be the ideal of the diagonal $\Delta$
of $C\x C$, and $p_1\:C\x C\to C$ be the first projection.  Then
$\cI_{\Delta}$ is a torsion-free rank-1 sheaf of degree $-1$ on $C\x
C/C$, and the tensor product $\cI_{\Delta}\ox p_1^*\cL$ defines
$A_\cL$.  Forming $A_\cL$ commutes with changing the base $S$,
and if the fibers of $C/S$ are not of arithmetic genus 0,
then $A_\cL$ is a closed embedding by \cite{AK80, (8.8), \p.108}.

Assume now that the geometric fibers of $C/S$ have only surficial
singularities (ones with embedding dimension 2), for example, double
points.  Then the projective $S$-scheme $\CJ^n$ is flat, and its
geometric fibers are integral local complete intersections; see
\cite{AIK76, (9), \p.8}.  Hence, the Picard scheme $\Pic_{\CJ^n/S}$
exists and is a disjoint union of quasi-projective $S$-schemes; see
Th\'eor\`eme~3.1, \p.232-06, in \cite{FGA}, and Corollary~(6.7)(ii),
\p.96, in \cite{AK80}.  So the Abel map induces an $S$-map,
        $$\textstyle A_\cL^*\:\Pic_{\CJ/S}\to \coprod_nJ^n.$$

As is customary \cite{FGA, \p.236-03}, let $\Pic^0_{\CJ/S}$ denote the
set-theoretic union of the connected components of the identity 0 in the
fibers of $\Pic_{\CJ/S}$, and let $\Pic^\tau_{\CJ/S}$ denote the set of
points of $\Pic_{\CJ/S}$ that have a multiple in $\Pic^0_{\CJ/S}$.  The
set $\Pic^\tau_{\CJ/S}$ is open; give it the induced scheme structure.

  The following theorem asserts that, if the geometric fibers of $C/S$
only have double points (of arbitrary order) as singularities, then
$\Pic^0_{\CJ/S}$ and $\Pic^\tau_{\CJ/S}$ are equal, and under
$A_\cL^*$, they are isomorphic to $J$.  This is our main result,
and its proof occupies the rest of the paper.

 \medbreak{\bf Theorem} \(Autoduality).  {\it
 Let $C/S$ be a flat projective family of integral curves.  Assume
its geometric fibers have double points at worst.
  Then\/ $\Pic^0_{\CJ/S}=\Pic^\tau_{\CJ/S}$.
  Furthermore,
 the Abel map induces an isomorphism,
        $$A_\cL^*\:\Pic^\tau_{\CJ/S}\risom J,$$
 which is independent of the choice of the invertible sheaf $\cL$ of
degree 1 on $C/S$; in fact, the isomorphism exists whether or not any
sheaf $\cL$ does.}

 \medbreak{\bf Proposition (\number\sctno.2)} \(Right inverse).
 \bgroup\it Let $C/S$ be a flat projective family of integral curves.
Assume its geometric fibers only have surficial singularities.
Then there exists a natural map,
        $$\beta\:J\to\Pic_{\CJ/S},$$
 whose formation commutes with base change, and whose image lies in the
subset $\Pic^0_{\CJ/S}$.  Furthermore, $A_\cL^*\circ\beta=
1_{J}$ for any $\cL$.
 \pf
 Choose an \'etale covering $S'/S$ such that the smooth locus of $C\x
S'/S'$ admits a section (such a covering exists by \cite{EGA, IV$_4$
17.16.3(ii), \p.106}).  Choose universal sheaves $\cI$ on $C\x\CJ\x
S'$ and $\cM$ on $C\x J\x S'$.  Form $C\x\CJ\x J\x S'$, and let
$p_{ijk}$ be the projection onto the product of the indicated factors.
Set
        $$\dM:=(\cD_{p_{234}}(p_{124}^*\cI\ox p_{134}^*\cM))^{-1}
        \ox\cD_{p_{234}}(p_{124}^*\cI)\on \CJ\x J\x S'$$
 where $\cD_{p_{234}}$ denotes the determinant of cohomology; see
 Section 6 in \cite{E99}, or \cite{KM76}.
 So $\dM$ is an invertible sheaf.  It defines the desired
map $\beta$ as we now prove.

The sheaf $\cI$ is determined up to tensor product with the pullback of
an invertible sheaf $\cN$ on $\CJ\x S'$.  So the projection formula
for the determinant of cohomology yields
  $$\eqalign{\cD_{p_{234}}(p_{124}^*\cI\ox p_{24}^*\cN\ox p_{134}^*\cM)
 &= \cD_{p_{234}}(p_{124}^*\cI\ox p_{134}^*\cM)\ox p_{13}^*\cN^{\ox m}\cr
        \cD_{p_{234}}(p_{124}^*\cI\ox p_{24}^*\cN)
 &= \cD_{p_{234}}(p_{124}^*\cI)\ox p_{13}^*\cN^{\ox n}\cr}$$
 where the $p$'s are the indicated projections and where $m$ and $n$ are
the Euler characteristics of $p_{124}^*\cI\ox p_{134}^*\cM$ and
$p_{124}^*\cI$ on the fibers of $p_{234}$ (thus $m$ and $n$ are locally
constant functions on $\CJ\x J\x S'$).  Now, $m=n$ because the
fibers of $p_{134}^*\cM$ have degree 0.  Therefore $\dM$ does not depend
on the choice of $\cI$.

Similarly, the sheaf $\cM$ is determined up to tensor product with the
pullback of an invertible sheaf $\cP$ on $J\x S'$.  Moreover, the
preceding argument shows that, if $\cM$ is replaced by its tensor product
with the pullback of $\cP$, then $\dM$ is replaced by its tensor product
with the pullback of $\cP^{\ox m}$.

Set $S'':=S'\x S'$.  There are two pullbacks of $\cI$ to $C\x\CJ\x
S''$, and both are universal sheaves.  Similarly, there are two
pullbacks of $\cM$ to $C\x J\x S''$, and both are universal sheaves.
Now, forming the determinant of cohomology commutes with changing the
base.  Therefore, by the preceding paragraphs, the two pullbacks of
$\dM$ to $\CJ\x J\x S''$ differ by tensor product with the
pullback of an invertible sheaf on $J\x S''$.  Hence $\dM$ defines a
map $\beta\:J\to \Pic_{\CJ/S}$.

Consider another choice of covering $S'_1/S$ and of sheaves $\cI_1$ and
$\cM_1$, and form the corresponding $\dM_1$.  Set $S'':=S'\x S'_1$.  Then
the pullbacks of $\cI_1$ and $\cI$ to $C\x\CJ\x S''$ are both
universal.  Similarly, the pullbacks of $\cM_1$ and $\cM$ to $C\x J\x
S''$ are both universal.  Hence, by the preceding argument, the
pullbacks of $\dM$ and $\dM_1$ to $\CJ\x J\x S''$ differ by
tensor product with the pullback of an invertible sheaf on $J\x S''$.
So $\dM_1$ and $\dM$ define the same map $\beta$.

Forming $\beta$ commutes with changing $S$ since forming the determinant
does.

The image of $\beta$ lies in $\Pic^0_{\CJ/S}$.  Indeed, we may change
the base to an arbitrary geometric point of $S$, and so work over an
algebraically closed field.  Then $J$ is integral.  So it suffices to
prove $\beta(0)=0$.  Now, we may choose $\cI$ on $C\x\CJ$ and $\cM$ on
$C\x J$.  Then the fiber $\cM(0)$ is equal to $\cO_C$.  Since forming
the determinant commutes with passing to the fiber, it follows that
$\dM(0)= \cO_{\CJ}$. So $\beta(0)=0$.

Finally, $A_\cL^*\circ\beta=1_{J}$.  Indeed, it suffices to check this
equation after changing the base to $S'$; so assume $S'=S$.  Then $\cI$
sits on $C\x\CJ$, and $\cM$ sits on $C\x J$.  So $A_\cL$ is defined by
$(1_C\x A_\cL)^*\cI$, as well as by $\cI_{\Delta}\ox p_1^*\cL$.  Hence
these two sheaves differ by tensor product with the pullback, along the
projection $p_2$, of an invertible sheaf on $C$.  It follows as above
from the properties of the determinant of cohomology that
 $$(A_\cL\x1_{J})^*\dM=
        (\cD_{p_{23}}(p_{12}^*\cI_{\Delta}\ox p_1^*\cL\ox p_{13}^*\cM))^{-1}
        \ox\cD_{p_{23}}(p_{12}^*\cI_{\Delta}\ox p_1^*\cL)\eqno\Cs2.1)$$
 on $C\x J$.  So both sides of this equation define the same map
$J\to\coprod_nJ^n$.

To evaluate the right-hand side of \Cs2.1), consider the natural sequence,
 $$0\to\cI_{\Delta}\to\cO_{C\x C}\to\cO_{\Delta}\to0.\eqno\Cs2.2)$$
 Pull it back to $C\x C\x J$, then tensor with $p_1^*\cL\ox
p_{13}^*\cM$ and with $p_1^*\cL$.   The additivity of the
determinant of cohomology now yields
 $$\eqalign{
 \cD_{p_{23}}(p_{12}^*\cI_{\Delta}\ox p_1^*\cL\ox p_{13}^*\cM)
 &=\cD_{p_{23}}(p_1^*\cL\ox p_{13}^*\cM)\ox(p_1^*\cL\ox\cM)^{-1},\cr
 \cD_{p_{23}}(p_{12}^* \cI_{\Delta}\ox p_1^*\cL)
 &=\cD_{p_{23}}(p_1^*\cL)\ox(p_1^*\cL)^{-1}
 .\cr}$$

Consider the following
 Cartesian square:
        $$\CD
                C\x J  @<p_{13}<<    C\x C\x J \\
              @VV p_2V @(\fiberbox) @VV p_{23}V \\
                   J    @<p_2<<        C\x J
        \endCD$$
 Forming the determinant commutes with changing the base.  So, on $C\x J$,
 $$\eqalign{
 \cD_{p_{23}}(p_1^*\cL\ox p_{13}^*\cM) &=p_2^*\cD_{p_2}(p_1^*\cL\ox \cM),\cr
 \cD_{p_{23}}(p_1^*\cL) &=p_2^*\cD_{p_2}(p_1^*\cL).\cr}$$
 Hence the right-hand side of \Cs2.1) differs from $\cM$ by
tensor product with the pullback of
an invertible sheaf on $J$.  Therefore,
$A_\cL^*\circ\beta=1_{J}$, and the proof is complete.

 \rmk3 In Proposition \Cs2), we made the hypothesis that the fibers of $C/S$
have surficial singularities, but we did not use the hypothesis directly
in the proof.  Rather, we used it indirectly to guarantee the existence
of the Picard scheme $\Pic_{\CJ/S}$.  Thus, the lemma is valid
whenever this Picard scheme exists, for example, when $S$ is the
spectrum of a field; see Corollaire 1.2 on \p.596 in \cite{R67}.

 \rmk4 Under the conditions of Proposition (2.2), assume that there is an
invertible sheaf $\cL$ of degree 1 on $C/S$.  Then the map
$\beta\:J\to\Pic_{\CJ/S}$ can be constructed in another and more
traditional way than that used in the proof of the lemma.  Namely,
$\beta$ can be constructed using the theta divisor associated to
$\cL$.  This is a divisor $\Theta_{\cL}$ on $\CJ$, and it may be
constructed as follows.

Use the notation of the proof of the proposition.  In addition, let $g$
denote the (locally constant) arithmetic genus of the fibers of $C/S$.
Now, on each fiber of the projection $p_{23}\: C\x\CJ\x S' \to \CJ\x
S'$, the restriction of $\cI\ox p_1^*\cL^{\ox g-1}$ has Euler
characteristic 0. It follows that, on $\CJ\x S'$, the invertible sheaf
        $$\cD_{p_{23}}(\cI\ox p_1^*\cL^{\ox g-1})$$
 has a canonical regular section; denote its divisor of zeros by
$\Theta'_{\cL}$.  Arguing as in the proof of the lemma, we can show that
$\Theta'_{\cL}$ descends to a divisor $\Theta_{\cL}$ on $\CJ$.

Let $\tau\: \CJ\x J \to \CJ$ be the multiplication map; it is defined by
$p_{124}^*\cI\ox p_{134}^*\cM$ on $C\x\CJ\x J\x S'$.  On $\CJ\x J\x S'$,
consider $\dM$, and on $\CJ\x J$, form the  sheaf,
    $$\cT:=\cO_{\CJ\x J}(p_1^*\Theta_{\cL}-\tau^*\Theta_{\cL}).$$
 We are about to construct a faithfully flat covering $S''/S'$ such that
the pullbacks of $\dM$ and $\cT$ are equal.  Each sheaf defines a map
from $J$ to $\Pic_{\CJ/S}$, and these two maps are equal after we change
the base to $S''$.  So, by descent theory, the two maps are equal to
begin with.  Therefore, since $\dM$ defines $\beta$, so does $\cT$.

To construct $S''/S'$, we may replace $S$ by $S'$, and so assume $S'=S$.
Moreover, we may assume that $S$ is affine, so Noetherian, and is
connected.
 After a further replacement of $S$, we may assume that the
smooth locus $\smC$ of $C/S$ admits a section $\sigma$; in fact, if we
replace $S$ by $\smC$, then the diagonal provides the desired section
$\sigma$.
  Fix $m$ so large that $\cL(m\sigma(S))$ is very ample.  Then,
again after replacing $S$, we can find a hyperplane section $H$ of $C$,
which is flat over $S$ and whose support lies in $\smC$.

Given any relative effective divisor $H_0$ on $\smC/S$ of relative
degree $n$, we can find a faithfully flat covering of $S$ such that,
after replacing $S$, we can find sections $\sigma_i$ of $\smC/S$ such
that
        $$H_0=\sigma_1(S)+\dots+\sigma_n(S).\eqno\Cs4.1)$$
 Indeed, $H_0/S$ is a faithfully flat covering, and after replacing $S$
by $H_0$, we have a canonical section $\sigma_1$ of $\smC/S$ whose image
is a subscheme of $H_0$ (in fact, $\sigma_1$ is simply the diagonal map
of the original $H_0/S$).  Form
        $$H_1:=H_0-\sigma_1(S).$$
 It is a relative effective divisor on $\smC/S$ of constant relative
degree $n-1$.  Hence, by induction, we may assume that, after replacing
$S$, we can find sections $\sigma_2,\dots,\sigma_n$ of $\smC/S$ such that
$H_1=\sigma_2(S)+\dots+\sigma_n(S)$. Then \Cs4.1) holds.
Taking $H_0$ to be $H$,
we conclude that we may assume that we have sections $\sigma_i$ of
$\smC/S$ such that
        $$\cL=\cO_C(T),\where
         T:=\sigma_1(S)+\dots+\sigma_{m+1}(S)-m\sigma(S).$$

Given a Cartier divisor $D$ on $C$, set $\cI(D):=\cI\ox p_1^*\cO_C(D)$
and
 $$\dM[D]:=\bigl(\cD_{p_{23}}(p_{12}^*\cI(D)\ox p_{13}^*\cM)\bigr)^{-1}
        \ox\cD_{p_{23}}(p_{12}^*\cI(D)).$$
 Then $\dM=\dM[0]$, and $\cT=\dM[(g-1)T]$.  Hence, it now suffices to
prove the following assertion: given any section $\rho$ of $\smC/S$, set
$R:=\rho(S)$ and let $E:=D+R$; then $\dM[D]$ and $\dM[E]$ differ by tensor
product with the pullback of an invertible sheaf on $J$.

To prove this assertion, consider the natural exact sequence,
        $$ 0 \to \cO_C(-R) \to \cO_C \to \cO_R \to 0.$$
 Pull it back to $C\x\CJ\x J$, then tensor with $p_{12}^*\cI(E)\ox
p_{13}^*\cM$ and with $p_{12}^*\cI(E)$.  Identify $R\x\CJ$ with
$\CJ$.  Additivity of the determinant of cohomology now yields
 $$\eqalign{ \cD_{p_{23}} (p_{12}^*\cI(E)\ox p_{13}^*\cM)=
        & \cD_{p_{23}}(p_{12}^*\cI(D)\ox p_{13}^*\cM)\ox
        p_1^*(\cI(E)|\CJ)\ox p_2^*(\cM|J)\cr
        \cD_{p_{23}}(p_{12}^*\cI(E))=
        & \cD_{p_{23}}(p_{12}^*\cI(D))\ox p_1^*(\cI(E)|\CJ).\cr} $$
 Hence $\dM[D]$ and $\dM[E]$ differ by tensor product with
$p_2^*(\cM|J)$.  So the assertion holds.  Thus $\cT$ defines
$\beta\:J\to\Pic_{\CJ/S}$.

\sct3 Theorem of the cube

\dsc1 Abel maps.  Let $C/S$ be a flat projective family of integral
curves, $m$ and $n$ integers.  The {\it Abel map\/} of bidegree
$(m,n)$ is defined to be the map,
        $$A_{C/S}\:\Hilb_{C/S}^m\x J^n_{C/S}\to \CJ^{n-m}_{C/S},$$
 given by tensoring the ideal of an $m$-cluster with a degree-$n$
invertible sheaf. We will often abbreviate $A_{C/S}$ by $A$.

More precisely, an $S$-map $t\:T\to \Hilb_{C/S}^m\x J^n$ corresponds to
a pair consisting of a flat closed subscheme $Y$ of $C\x T/T$ with
length-$m$ fibers and of an invertible sheaf $\cL'$ on $C\x T'/T'$ with
degree-$n$ fibers, where $T'/T$ is an \'etale covering, such that the
two pullbacks of $\cL'$ to $C\x T''$ are equal, where $T''\big/\,
T'\x_TT'$ is an \'etale covering.  Let $\cI'$ denote the ideal of $Y\x
T'$.  Then $\cI'\ox\cL'$ is a torsion-free rank-1 sheaf of degree $n-m$
on $C\x T'/T'$, and its two pullbacks to $C\x T''$ are equal.  Hence
$\cI'\ox\cL'$ defines a map $A(t)\:T\to\CJ^{n-m}$.

The Abel map is smooth when the geometric fibers of $C/S$ have double
points at worst, thanks to the following more general fact.
\smallskip\hang
 {\smc (Smoothness)\enspace \it If all the fibers of $C/S$ are Gorenstein,
then the Abel map $A$ is smooth.}
 \smallskip\noindent
 This fact is proved in Corollary (2.6) of \cite{EGK}, as an application
of an even more general statement, Theorem (2.4) of
 \cite{EGK}.

 \lem2 Let $C/S$ be a flat projective family of integral curves.  Assume
there is a universal sheaf $\cI$ on $C\x\CJ^1$.  Set $P:=\IP(\cI)$.  Let
$Z\subset P$ be the preimage of $C\x J^1$.  Then the structure map of
$P$ induces an isomorphism $Z\risom C\x J^1$, and there is a map
$\zeta\:P\to\CJ$ extending the Abel map $A\:C\x J^1\to \CJ$.
 \pf Let $\rho\:P\to C\x\CJ^1$ be the structure map; say
$\rho=(\rho_1,\rho_2)$.  Set $\gamma:=(\rho_1,1_P)$ and $\theta:=
1_C\x\rho_2$.  Then $\rho$ factors as follows:
        $$\rho\:P \TO\gamma C\x P \TO\theta C\x\CJ^1.$$
 So there is a natural isomorphism $\gamma^*\theta^*\cI\risom\rho^*\cI$.
Let $q\:\theta^*\cI\to\gamma_*\rho^*\cI$ be its adjoint.  In other
words, $\gamma$ is the graph map of $\rho_1$, and its image, $Y$ say, is
the graph subscheme; in these terms, $q$ is equal to the natural
quotient map $\theta^*\cI\onto\theta^*\cI|Y$.

Let
 $u\:\rho^*\cI\to\cO_P(1)$
 be the universal map, and form the composition,
   $$r\:\theta^*\cI \TO q \gamma_*\rho^*\cI
         \TO{\gamma_*u} \gamma_*\cO_P(1).$$
 Then $r$ is a surjective map between $P$-flat sheaves on $C\x P$.
Set $\cJ:=\cKer(r)$.  Then $\cJ$ is flat too, and forming it commutes
with passing to the fibers.  Hence $\cJ$ is a torsion-free rank-1 sheaf of
degree 0 on $C\x P$.  So it defines a map $\zeta\:P\to\CJ$.

Since $\cI|C\x J^1$ is invertible, $\rho$ restricts to an isomorphism
$Z\risom C\x J^1$, and $u\:\rho^*\cI\to\cO_P(1)$ restricts to an
isomorphism.  Hence the exact sequence,
        $$0\To\cJ\To\theta^*\cI\TO r\gamma_*\cO_P(1)\To0,$$
 is equal, on $C\x Z$, to the tensor product of $\theta^*\cI$ with the
basic sequence,
        $$0\To\cI_Y\To\cO_{C\x P}\To\cO_Y\To 0,$$
 where $\cI_Y$ is the ideal of the graph subscheme $Y$.  However, this
sequence is the pullback under $1_C\x\rho_1\:C\x P\to C\x C$ of the
basic sequence of the diagonal, (2.2.2).  Hence $\cJ|C\x Z$ defines
$A\circ(\rho|Z)$.  Thus $\zeta$ extends $A\:C\x J^1\to \CJ$,
and the proof is complete.

 \rmk3 In the proof of Proposition~\Cs2), consider the map $r$.  Being
surjective, $r$ defines a map,
        $$P\to\Quot^1_{\cI/C\x\CJ^1/\CJ^1},$$
 and it is not hard to see that this map is an isomorphism.  (In fact,
there is nothing special about $\cI/C\x\CJ^1/\CJ^1$.  This is a simple
general phenomenon.  See \cite{K87, (2.2), \p.109}.)  Thus
$\zeta\:P\to\CJ$ is a universal lifting of $A\:C\x J^1\to
\CJ$ in the following sense: given a map $T\to\CJ$ defined by a
degree-0 torsion-free rank-1 subsheaf of a degree-1 torsion-free rank-1
sheaf on $C\x T$, there exists a unique map $T\to P$ such that the pair
of sheaves is the pullback of the pair consisting of $\cJ$ and
$\theta^*\cI$.

 \lem4 Let $C/S$ be a flat projective family of integral curves, and
$T\to C$ an $S$-map.  Let $\Gamma$ be the graph subscheme of $C\x T$,
and $\cI_\Gamma$ its ideal.  Set $W:=\IP(\cI_\Gamma)$, and let
$\psi\:W\to C\x T$ be the structure map.  Assume that the geometric
fibers of $C/S$ only have surficial singularities.  Then $W/T$ is flat,
and
 $$\psi_*\cO_W=\cO_{C\x T}\and R^i\psi_*\cO_W=0\for i\ge1.\eqno\Cs4.1)$$
 \pf
 Without loss of generality, we may replace $S$ by $\Spec\cO_{S,s}$
where $s$ is an arbitrary point of $S$.
By \cite{EGA, O$_{\rm III}$ 10.3.1,
\p.20}, there exists a flat local $\cO_{S,s}$-algebra whose residue
field is any given extension of the field of $s$, and we may replace $S$
by the spectrum of this algebra.  Thus we may assume that $S$ is a local
scheme with closed point $s$ whose residue field $k(s)$ is algebraically
closed.

Embed $C/S$ in a projective space $\IP^N_S$ for some $N$.  Let $\cH$ be
the ideal of $C$.  Since $C/S$ is flat, $\cH(s)$ is the ideal of $C(s)$,
and $\cH$ is flat.  Also, for $m\gg0$, the base change map is an
isomorphism,
        $$H^0(\cH(m))\ox k(s)\risom H^0(\cH(s)(m)).$$

Fix $m\gg0$, and take $N-2$ general sections of $\cH(s)(m)$.  Via the
above isomorphism, lift the sections back to sections of $\cH(m)$, and
form the scheme $F$ of common zeros of the lifts.  Then $F\supset C$.
Moreover, increasing $m$ if necessary, we may assume that $F(s)$ is a
smooth surface since every singularity of $C$ is surficial
(see the proofs of (7)--(9) in \cite{AK79} for example).
Hence $F/S$ is a smooth family of surfaces.

 Consider the nested sequence of subschemes,
        $$\Gamma\subset C\x T \subset F\x T.$$
 Since $F/S$ is smooth and $\Gamma$ is a graph, $\Gamma$ is regularly
embedded in $F\x T$, say with ideal $\cJ$.  Hence the symmetric algebra
of $\cJ$ is equal to its Rees algebra by Micali's theorem \cite{M60,
\p.1955}.  So $\IP(\cJ)$ is equal to the blowup $B$ of $F\x T$ along
$\Gamma$.  Since $F\x T/T$ and $\Gamma/T$ are both smooth, so is $B/T$.

Denote the sheaf of ideals of $C\x T$ in $F\x T$ by $\cK$, and form the
exact sequence,
        $$0\to\cK\to\cJ\to\cI_\Gamma\to0.$$
 It gives rise to the following exact sequence of sheaves of graded modules
over the symmetric algebra $\cSym(\cJ)$ (see \cite{E95, \p.571}):
        $$0\to\cK\cdot\cSym(\cJ)[-1]\to\cSym(\cJ)
        \to\cSym(\cI_\Gamma)\to0.$$
 Taking associated sheaves yields the following exact sequence on $B$:
        $$0\to\cK\cdot\cO_B(-1)\to\cO_B\to\cO_W\to0.\eqno\Cs4.2)$$
Thus $\cK\cdot\cO_B(-1)$ is the ideal of $W$ on $B$.

Since $F/S$ is smooth and $C/S$ is flat, $C$ is a Cartier divisor on
$F$.  Hence $\cK$ is invertible.  Thus $W$ is a Cartier divisor on $B$;
in fact, $W=D-E$ where $D$ is the preimage of $C\x T$ in $B$ and where
$E$ is the exceptional divisor.  Since $W$ remains a Cartier divisor on
the fibers of $B/T$ and since $B/T$ is flat, $W/T$ is flat.

Consider the blowup map $\beta\:B\to F\x T$.  Since $\cK$ is invertible,
the projection formula yields, for every $i$,
        $$R^i\beta_*\cK\cdot\cO_B(-1)=\cK\cdot R^i\beta_*\cO_B(-1).$$
 Hence Assertion \Cs4.1) follows from the long exact sequence of higher
direct images of $\beta$ associated to the sequence \Cs4.2) and from the
following formulas:
        $$\beta_*\cO_B(n)=\cJ^n\and R^i\beta_*\cO_B(n)=0
         \for i\ge1 \and n\ge-1,\eqno\Cs4.3)$$
 where $\cJ^n=\cO_{F\x T}$ for $n\le0$ by convention.  These formulas
are proved next (and the proof applies more generally to any blowup
along a regularly embedded center of codimension at least 2).

Since $B=\IP(\cJ)$, restricting the base yields $E=\IP(\cJ/\cJ^2)$.  Now,
$\cJ/\cJ^2$ is locally free.  Hence $R^i\beta_*\cO_E(n)=0$ for
$i\ge1$ and $n\ge-1$ by
Serre's computation.  So the long exact sequence associated to the
sequence,
        $$0\to\cO_B(n+1)\to\cO_B(n)\to\cO_E(n)\to0,\eqno\Cs4.4)$$
 yields a surjection $R^i\beta_*\cO_B(n+1)\onto R^i\beta_*\cO_B(n)$ for
$i\ge1$ and for $n\ge-1$.  By Serre's theorem, $R^i\beta_*\cO_B(n)$
vanishes for $i\ge1$ and $n\gg0$.  Hence, by descending induction on $n$, it
vanishes for $i\ge1$ and $n\ge-1$.

Sequence \Cs4.4) also gives rise to the following commutative diagram:
        $$\CD
  0@>>>   \cJ^{n+1}     @>>>   \cJ^n       @>>> \cJ^n/\cJ^{n+1} @>>>0\\
  @.      @VVV                 @VVV                     @VVV \\
  0@>>>\beta_*\cO_B(n+1)        @>>> \beta_*\cO_B(n)@>>> \beta_*\cO_E(n)
        \endCD$$
 Since $E=\IP(\cJ/\cJ^2)$, the right vertical map is an isomorphism for
$n\ge-1$ by Serre's computation again.  The left vertical map is an
isomorphism for $n\gg0$ by Serre's theorem.  Hence, by descending
induction on $n$, the middle vertical map is an isomorphism for
$n\ge-1$.  Thus Formulas \Cs4.3) hold, and the proof is complete.

 \medbreak{\bf Lemma (\number\sctno.5)} \(Generalized theorem of the
cube).  \bgroup\it Let $S$ be a connected and locally Noetherian scheme.
Let $g\:Y\to S$ and $f\:X\to Y$ be flat and proper maps.  Let
$\sigma\:S\to Y$ and $\tau\:Y\to X$ be sections of $g$ and $f$.
Assume\smallbreak
 \item i that  $\cO_S=g_*\cO_Y$  and $\cO_Y=f_*\cO_X$ hold universally, and
 \item ii that, for every closed point $s\in S$,
the natural map on the fibers is
injective:
        $$w\:H^0(Y(s), R^1f(s)_*\cO_{X(s)}) \into
        H^1(f(s)^{-1}\sigma(s),\cO_{f(s)^{-1}\sigma(s)}).$$
 \smallbreak\noindent
 Let $s_0\in S$ and $\cL$ be an invertible sheaf on $X$.  If
the three restrictions,
       $$\cL|X(s_0),\ \cL|\tau(Y),\and\cL|f^{-1}\sigma(S),$$
 are trivial, then $\cL$ is trivial.
 \pf
 It is not hard to adapt, mutatis mutandis, Mumford's proof of his
similar theorem on \p.91 in \cite{M70}.  It is a straightforward job,
except at the beginning and at the end.  At the beginning, Mumford uses
his proposition on \p.89 to obtain the existence of a maximum closed
subscheme $T$ of $S$ carrying an invertible sheaf $\cN$ such that
$h^*\cN=\cL|h^{-1}T$ where $h:=gf$.  (In fact, forming $T$ commutes with
base-changing $S$.)  Mumford's construction of $T$ is not hard to adapt;
here's the idea.

Since $h$ is flat and proper, by \cite{EGA, III$_2$ 7.7.6, \p.69},
there are coherent sheaves $\cM$ and $\cN$ on $S$ such that, for
every coherent sheaf $\cF$ on $S$, we have
        $$h_*(\cL\ox h^*\cF)=\cHom(\cM,\cF)
        \and h_*(\cL^{-1}\ox h^*\cF)=\cHom(\cN,\cF).$$
 Take $T:=\Supp(\cM)\bigcap\Supp(\cN)$ using the annihilators of $\cM$
and $\cN$ to define the scheme structures on their supports.

Over a point $t\in T$, the restrictions $\cL|X(t)$ and $\cL^{-1}|X(t)$
each have a nonzero section.  Compose the first section with the dual of
the second, obtaining a nonzero map $\cO_{X(t)}\to\cO_{X(t)}$.  This map
is given by multiplication by a scalar because $h_*\cO_X=\cO_T$ holds
universally.  Hence $\cL|X(t)$ is trivial.  Therefore, on a neighborhood
of $t$, one element suffices to generate $\cN$.  It follows that $\cN|T$
is an invertible $\cO_T$-module and that $h^*\cN|h^{-1}T=\cL|h^{-1}T$.

At the end of the proof of his theorem, Mumford uses the K\"unneth
formula.  Instead, we must use a related injection, which we get
as follows.  For each closed point $s\in S$,
form the exact sequence of terms of low
degree of the Leray spectral sequence:
 $$H^1(Y(s),f(s)_*\cO_{X(s)}) \TO u H^1(X(s),\cO_{X(s)}) \TO v
        H^0(Y(s),R^1f(s)_*\cO_{X(s)}). $$

Since $f(s)_*\cO_{X(s)}=\cO_{Y(s)}$ by Assumption (i), the section
$\tau(s)\:Y(s)\to X(s)$ of $f(s)$ yields a map,
        $$u'\: H^1(X(s),\cO_{X(s)}) \to H^1(Y(s),\cO_{Y(s)}),$$
 splitting $u$.  Hence, by Assumption (ii), the map $(u',w\circ v)$ is
an injection,
  $$H^1(X(s),\cO_{X(s)}) \into H^1(Y(s),\cO_{Y(s)})\textstyle
   \bigoplus H^1(f(s)^{-1}\sigma(s),\cO_{f(s)^{-1}\sigma(s)}).$$
 This injection works in place of the K\"unneth formula.

 \lem6 Let $C/S$ be a flat projective family of integral curves with
surficial singularities, and let $\cP$ be an invertible sheaf on $\CJ$.
Assume that $S$ is connected, and that $S$ contains a point $s_0$ such
that the fiber $\cP(s_0)$ is trivial.  Form the Abel pullback
$A^*\cP$ on $C\x J^1$.  Then the following two assertions
hold.
 \part 1 The pullback $A^*\cP$ defines a map $a\:J^1\to J$,
and $a$ factors uniquely through the structure map $J^1\to S$.
 \part 2 If the smooth locus $\smC$ of $C/S$ admits a section, then
there are invertible sheaves $\cM_1$ on $C$ and $\cM_2$ on $J^1$ such
that
        $$A^*\cP=p_1^*\cM_1\ox p_2^*\cM_2,$$
 where the $p_i$ are the projections.
 \pf
 Consider Part (1).  A priori, $A^*\cP$ defines a map from
$J^1$ to $\Pic_{C/S}$.  However, the image lies in the open subscheme
$J$ because $A^*\cP(s_0)$ is trivial and $S$ is connected.

If there is an $S$-factorization $a\:J^1\TO bS\TO cJ$, then $b$ must
be the structure map.  So $b$ is faithfully flat.  Hence, by descent
theory, $c$ is uniquely determined.

Assume for a moment that the smooth locus $\smC$ of $C/S$ admits a
section, and that Part (2) holds.  Then $\cM_1$ defines an $S$-map
$c\:S\to J$ such that $a=cb$, where $b$ is the structure map.

In general, there is an \'etale covering $S'/S$ such that $\smC\x
S'/S'$ admits a section.  So, by the reasoning above, $a\x S'$ factors
through a unique map $c'\:S'\to J\x S'$.  Set $S'':=S'\x S'$.  Then $a\x
S''$ factors through both $c'\x S'$ and $S'\x c'$.  So the latter two
maps are equal.  Hence $c'$ descends to a suitable map $c\:S\to J$.
Thus Part (1) follows from Part (2).

To prove Part (2), assume from now on that $\smC$ admits a section.
Then there exists a universal sheaf $\cI$ on $C\x\CJ^1$.  Set
$P:=\IP(\cI)$.  Let $\rho\:P\to C\x\CJ^1$ denote the structure map, and
        $$\rho_1\:P\to C\and\rho_2\:P\to\CJ^1$$
 the natural maps.  Now, $A\:C\x J^2\to \CJ^1$ is smooth; see
(3.1).  So its image $V$ is open.  Also $V\supset J^1$.  Set
        $$X:=\rho^{-1}(C\x V)\and Z:=\rho^{-1}(C\x J^1).$$
 Then $\rho$ induces an isomorphism $Z\risom C\x J^1$; so the $\rho_i$
extend the projections $p_i$.

Say $\theta\:S\to C$ is the section, let $Q$ be its image, and set
$\cL:=\cO_C(Q)$.  Then $\cO_C(Q)$ defines a section $\phi\:S\to \CJ^1$,
whose image lies in $J^1$.  So $\phi$ defines a section $\xi_1\:C\to
P$ of $\rho_1$ because $\rho$ is an isomorphism over $C\x J^1$.
Moreover, $\theta$ defines a section $\xi_2\:\CJ^1\to P$ of $\rho_2$
because $\rho$ is an isomorphism over $\smC\x\CJ^1$.

Consider the map $\zeta\:P\to\CJ$ of \Cs2) extending
$A\:C\x J^1\to\CJ$.  Set
        $$\vcenter{\openup1\jot
        \halign{$\hfil#$&${}#\hfil$&$#$\hfil\cr
        \cQ_2&:=\zeta^*\cP,   &\cN_2:=\xi_2^*\cQ_2;\cr
 \cQ_1&:=\cQ_2\ox \rho_2^*\cN_2^{-1},&\cN_1:=\xi_1^*\cQ_1;\cr
 \cQ_0&:=\cQ_1\ox \rho_1^*\cN_1^{-1}=\cQ_2\ox\rho_2^*&\cN_2^{-1}
        \ox \rho_1^*\cN_1^{-1}.\cr}}$$
 Notice that, by hypothesis and by construction, the restrictions,
        $$\cQ_0|P(s_0)\and\cQ_0|\xi_2(\CJ^1)
                \and\cQ_0|\rho_2^{-1}\phi(S),\eqno\Cs6.1)$$
 are trivial.  It suffices to prove that $\cQ_0|X$ is trivial.

Set $q:=\rho_2|X$, so $q\:X\to V$.  Then $q$ is flat.  Indeed, let
$\Delta$ be the diagonal subscheme of $C\x C$, and $\cI_\Delta$ its
ideal.  Set $W:=\IP(\cI_\Delta)$.  Then there is a
 Cartesian square
        $$\CD
         X              @<<<      W\x J^2 \\
        @VVV        @(\fiberbox)    @VVV \\
        C\x V @<1\x A<< C\x C\x J^2
        \endCD$$
 because, on $C\x C\x J^2$, the pullback of $\cI|C\x V$ is equal to the
tensor product of the pullback of $\cI_\Delta$ with an invertible sheaf
(namely, the pullback of a universal sheaf on $C\x J^2$).  So there is a
 Cartesian square
        $$\CD
         X      @<<<       W\x J^2 \\
        @VVqV @(\fiberbox)  @VVV \\
        V @<A<< C\x J^2
        \endCD$$
 The right vertical map is flat thanks to Lemma \Cs4) applied with $C$
for $T$.  Also, $A\:C\x J^2\to V$ is faithfully flat, being smooth and
surjective.  Hence $q$ is flat.

Similarly, $q_*\cO_X=\cO_V$ holds universally.  Indeed, by the preceding
argument, this statement results from the corresponding statement about
the composition $W\to C\x C\to C$.  The corresponding statement results
from Lemma \Cs4) and a basic fact: $p_*\cO_C=\cO_S$ holds universally
where $p\:C\to S$ is the structure map.

The generalized theorem of the cube, Lemma \Cs5), does not apply to the
triple $X/V/S$ because $V/S$ is not proper.  However, $V$ is swept out
by copies of $C$ because $A\: C\x J^2\to V$ is surjective.
So we can circumvent the obstacle as follows.

Since $q\:X\to V$ is flat and proper, and since $q_*\cO_X=\cO_V$ holds
universally, there exists a maximum closed subscheme $Y$ of $V$ carrying
an invertible sheaf $\cH$ such that $q^*\cH=\cQ_0|q^{-1}Y$, and forming
$Y$ commutes with base-changing $V$; see the proof of Lemma~\Cs5).  In
fact, $\cH$ is trivial because $\cQ_0|\xi_2(V)$ is trivial.  Hence
$\cQ_0|X$ is trivial if $Y=V$.  Since forming $Y$ commutes with
base-changing $V$, it suffices to construct a faithfully flat map $U\to
V$ such that the pullback of $\cQ_0$ to $X\x_VU$ is trivial.

Given $n\ge0$, let $C_n$ denote the $n$-fold self-product of $\smC$.
Let $\gamma_n\:C_n\to J^2$ be the map that sends $n$ $T$-points of
$\smC$, say with graph images $Q_1,\dots,Q_n$, to the $T$-point of $J^2$
representing the following invertible sheaf on $C\x T$:
        $$\cO((n+2)Q\x T-Q_1-\cdots-Q_n).$$
 Finally, set $\delta_n:=A\circ(1_C\x\gamma_n)$,
so $\delta_n\:C\x C_n\to V$.

Consider the factorization $\gamma_n\:C_n\TO b \Hilb^n_{\smC/S}\TO cJ^2$.
The map $b$ is faithfully flat for every $n$, and $c$ is faithfully flat
for $n\gg0$.  So $\gamma_n$ is faithfully flat for $n\gg0$.  Hence,
$\delta_n$ is faithfully flat for $n\gg0$, since
$A\: C\x J^2\to V$ is faithfully flat.

Thus, it suffices to prove that the pullback of $\cQ_0$ to $X\x_V(C\x
C_n)$ is trivial for every $n$.  To do so, we apply Lemma \Cs5) to
$X\x_V(C\x C_n)\big/(C\x C_n)\big/C_n$ with the sections,
        $$\sigma:=\theta\x1_{C_n} \and \tau:=\xi_2\x1_{C\x C_n}.$$
 Assumption (i) of \Cs5) holds as $q_*\cO_X=\cO_V$ and $p_*\cO_C=\cO_S$
hold universally.

To check Assumption (ii), let $t\in C_n$.  Set $T:=C\ox k(t)$ and
$W_T:=W\ox k(t)$.  Denote by $\psi\:W_T\to C\x T$ the structure map, and
by $\psi_2\:C\x T\to T$ the second projection.  Then Lemma~\Cs4) implies
that $\psi_*\cO_{W_T}=\cO_{C\x T}$ and $R^1\psi_*\cO_{W_T}=0$.  Hence,
using the Leray spectral sequence and making the substitution, we get
        $$R^1(\psi_2\psi)_*\cO_{W_T}=R^1\psi_{2*}(\psi_*\cO_{W_T})
        =R^1\psi_{2*}\cO_{C\x T}.$$
 Now, $C\x T=T\x T$.  So, commuting cohomology with flat base change, we
get
        $$R^1\psi_{2*}\cO_{C\x T}=H^1(T,\cO_T)\ox\cO_T.$$
 Therefore, the following natural map is an isomorphism:
        $$H^0(T,\,R^1(\psi_2\psi)_*\cO_{W_T})\risom H^1(T,\cO_T).$$
 It now follows immediately that Assumption (ii) holds.

It remains to check the triviality of the three appropriate pullbacks of
$\cQ_0$.  First, $C_n$ is connected and maps onto $S$ because the fibers
of $C/S$ are geometrically integral.  Moreover,  $\cQ_0|P(s_0)$ is
trivial by \Cs6.1).  Hence, $C_n$ contains a point $t_0$ that maps to
$s_0$, and the pullback of $\cQ_0$ to $X\ox k(t_0)$ is trivial.

Second, consider the pullback of $\cQ_0$ to $\tau(C\x C_n)$.  This
pullback is trivial because $\cQ_0|\xi_2(\CJ^1)$ is trivial by \Cs6.1).

Finally, consider the pullback of $\cQ_0$ to $X\x_V \sigma(C_n)$.  To
begin, suppose $n=0$.  Now, $C_0=S$.  So $\sigma=\theta$ and
$\delta_0\theta=\phi$.  Hence $X\x_V \sigma(C_0)$ is equal to
$\rho_2^{-1}\phi(S)$.  However, the pullback of $\cQ_0$ to the latter
scheme is trivial by \Cs6.1).

Proceeding by induction on $n$, suppose that the pullback of $\cQ_0$ to
$X\x_V \sigma(C_n)$ is trivial.  Then Lemma~\Cs5) implies that the
pullback of $\cQ_0$ to $X\x_V(C\x C_n)$ is trivial.  Hence so is the
pullback to $X\x_V(\smC\x C_n)$.  However, the latter scheme is equal to
$X\x_V \sigma(C_{n+1})$, as is easy to see.  The proof of the lemma is
now complete.

 \prp7 Let $C/S$ be a flat projective family of integral curves with
surficial singularities.  Assume $S$ is connected, and let $U$ be the
connected component of $\Pic_{\CJ/S}$ containing the zero section.
Then there exists a natural map $c\:U\to J$ such that
$c\circ\beta=1_{J}$, where $\beta$ is the map of Proposition~\(2.2).
Furthermore, given any invertible sheaf $\cL$ of degree $1$ on $C/S$,
the map $A_\cL^*$ of Subsection~\(2.1) restricts to $c$; in
particular, $A_\cL^*|U$ is independent of $\cL$.
 \pf
 The structure sheaf $\cO_C$ defines a section of $\CJ/S$.  Hence,
$\CJ\x U$ admits a universal invertible sheaf
 $\cP$; see Prop.\ts2.1 on Page 232-04 of \cite{FGA}.  Also,
 for every point $u_0$ of $U$ on the identity section,
$\cP(u_0)$ is trivial.  On $C\x J^1\x U$, form $(A\x1_U)^*\cP$.  This
pullback defines a map $a\:J^1\x U\to J$.  It factors through a map
$c\:U\to J$ by virtue of Part (1) of Lemma~\Cs6) applied to $C\x U/U$
and $\cP$.

Given $\cL$, let $[\cL]\in J^1(S)$ represent it.  Then the fiber
$a([\cL])\:U\to J$ is equal to $c$ on the one hand, and to
$A_\cL^*|U$ on the other.  Thus the last assertion holds.

The equation $c\circ\beta=1_{J}$ follows.  Indeed, it suffices to
check this equation after making an \'etale base change $S'/S$.  After a
suitable such base change, there exists a $\cL$.  Then, by what we just
proved, $c$ is equal to $A_\cL^*$, and so Proposition~(2.2) yields the
asserted equation.  The proof is now complete.

\sct4 Proof and extension

 \dsc1 Proof of the autoduality theorem of {\(2.1)}.  For a moment, make
the following assumption: for each geometric point $s$ of $S$ and for
some invertible sheaf $\cL_s$ of degree 1 on the curve $C(s)$, the Abel
map induces an isomorphism,
        $$A_{\cL_s}^*\:\Pic^\tau_{\CJ(s)} \risom J(s).$$
 This case of the theorem implies the general case, as we'll now prove.

Set $U:=\Pic^\tau_{\CJ/S}$ and consider the map
$\beta\:J\to\Pic_{\CJ/S}$ of Proposition~(2.2). The image of $\beta$
lies in $\Pic^0_{\CJ/S}$, hence also in $U$.  Since
$A_\cL^*\circ\beta= 1_{J}$ for any invertible
sheaf $\cL$ of degree 1 on $C/S$, we have to prove that
$\beta\:J\to U$ is an isomorphism.  To do so, we may change the base to
an \'etale covering of $S$.  Thus we may assume that the smooth locus of
$C/S$ admits a section.  Then $C/S$ does carry an invertible sheaf $\cL$
of degree 1.  Hence $\beta\:J\to U$ is a right inverse, so a closed
embedding.  Since $J$ is flat over $S$, it follows that $\beta\:J\to U$ is
an isomorphism, being one on each geometric fiber.  Therefore,
$U=\Pic^0_{\CJ/S}$.  Thus to prove the theorem, we may
assume that $S$ is the spectrum of an algebraically closed field.

Proceed by induction on the difference $\delta$ between the arithmetic
genus and the geometric genus of $C$.  First, assume $\delta=0$.  Then
$C$ is smooth.  Hence $J$ is complete, so an Abelian variety.  Given any
Abelian variety $G$, in the theorem on \p.125 of \cite{M70}, Mumford
proves that the scheme $\Pic^0_G$ is a quotient of $G$ by a finite
group; hence, $\Pic^0_G$ is integral and has the same dimension as $G$.
Moreover, $\Pic^0_G$ is equal to $\Pic^\tau_G$ by Corollary 2 on \p.178
of \cite{M70}.  Now, $A_\cL^*\circ\beta= 1_{J}$ for any invertible
sheaf $\cL$ of degree 1 on $C$ by Proposition~(2.2).  So $\beta$ is a closed
embedding of $J$ in $\Pic^0_{J}$.  Hence $\beta$ is an isomorphism.
Thus the theorem holds when $\delta=0$.

Assume $\delta\ge1$ from now on.  Fix an invertible sheaf $\cL$ of
degree 1 on $C$.  Then $A_\cL^*\circ\beta= 1_{J}$ by Proposition~(2.2).
So $A_\cL^*$ is an epimorphism.  We must prove it is a
monomorphism.  So let $\phi\:T\to U$ be an $S$-map of finite type, say
arising from the invertible sheaf $\cN$ on $\CJ\x T$.  Assume that
$(A_\cL\x1_T)^*\cN$ is equal to the pullback of an invertible sheaf
on $T$.  We must prove that $\cN$ is the pullback of an invertible sheaf
on $T$.  To do so, we may assume that $T$ is connected.

First set $U^0:=\Pic^0_\CJ$ and assume $\phi(T)\subset U^0$.  Note that
$U^0$ is an open subscheme of $U$ because we are now working over an
algebraically closed field.  Let $\cP$ be a universal invertible sheaf
on $\CJ\x U^0$.  Then there is an invertible sheaf $\cT$ on $T$ such
that
        $$(1\x\phi)^*\cP=\cN\ox q_2^*\cT \on \CJ\x T,$$
 where $q_2\:\CJ\x T\to T$ is the projection.

Let $u_0\in U^0$ denote the identity; so $\cP(u_0)$ is trivial.  Hence
Lemma~(3.6) applies to $C\x U^0/U^0$ and $\cP$.  Part (2) of the lemma
implies that
 $(A\x1)^*\cP$ is equal
on $C\x J^1\x U^0$ to the
tensor product of the pullbacks of invertible sheaves on $C\x U^0$ and
$J^1\x U^0$.  Therefore, there are invertible sheaves $\cN_1$ on $J^1\x T$
and $\cN_2$ on $C\x T$ such that
    $$(A\x1)^*\cN=q_{23}^*\cN_1\ox q_{13}^*\cN_2\on C\x J^1\x T,$$
 where $q_{23}$ and $q_{13}$ are the projections.

  Let $[\cL]\in J^1$ represent $\cL$.  Then the equation above yields
  $$(A_\cL\x1)^*\cN=((A\x1)^*\cN)[\cL]
        =(q_2^*\cN_1[\cL])\ox\cN_2\on C\x T.$$
 By assumption, the term on the left is the pullback of an invertible
sheaf on $T$; hence, so is $\cN_2$.  Therefore,
        $$(A\x1)^*\cN=q_{23}^*\cR\eqno\Cs1.1)$$
 for some invertible sheaf $\cR$ on $J^1\x T$.

Since $\delta\ge1$, there is a double point $Q\in C$.  Let $\vf\:\dr
C\to C$ be the blowup at $Q$.  Then there is a natural scheme $P$, known
as the {\it presentation scheme,} and there are natural maps
$\kappa\:P\to \CJ_C$ and $\pi\:P\to \CJ_{\dr C}$; see
Subsections~(3.1)~and~(3.3), in \cite{EGK}, or \cite{AK90}.  Since $Q$
is a double point, $P$ is a $\IP^1$-bundle over $\CJ_{\dr C}$; see
Theorem~(6.3) in \cite{EGK}.  Now, since $\cN$ arises from a map $T\to
U^0$, the pullback $(\kappa\x1)^*\cN$ restricts on each $\IP^1$ to a
sheaf of degree 0, hence is the pullback of an invertible sheaf $\dr\cR$
on $\CJ_{\dr C}\x T$.

Set $\dr\cL:=\vf^*\cL$.  By Lemma~(6.4) in \cite{EGK} the singularities
of $\dr C$ are only double points.  Hence, by Corollary~(5.5) in \cite{EGK},
there is a map $\Lambda\:\dr C\x J^1_C\to P$ making the following two diagrams
commute:
 $$\CD
        \dr C\x J_C^1       @>\Lambda>>             P       \\
           @V 1\x\vf^*VV                         @V\pi VV        \\
         \dr C\x J_{\dr C}^1  @>A_{\dr C}>> \CJ_{\dr C}
  \endCD\qquad\CD
        \dr C\x J_C^1  @>\Lambda>>          P       \\
           @V\vf\x1VV                    @V\kappa VV     \\
          C\x J_C^1    @>A_C>> \CJ_C
  \endCD\eqno\Cs1.2)$$
 Those diagrams yield these equations:
 $$(A_{\dr\cL}\x1)^*\dr\cR=(\Lambda_\cL\x1)^*(\kappa\x1)^*\cN
        =(\vf\x1)^*(A_\cL\x1)^*\cN=p_2^*\cR[\cL],$$
 where $p_2\:\dr C\x T\to T$ is the projection, and $\Lambda_\cL$ is the
composition of $\Lambda$ with the map $\dr C\to\dr C\x J^1_C$ defined by
$\cL$.

 By induction, autoduality holds for $\dr C$.  Hence $\dr\cR$ is the
pullback of an invertible sheaf on $T$, whence so is $(\kappa\x1)^*\cN$.
Invert the latter sheaf on $T$, pull it back to $\CJ_C\x T$, tensor with
$\cN$, and use the product to replace $\cN$.  Thus we may assume that
$(\kappa\x1)^*\cN$ is trivial.  So $(\Lambda\x1)^*(\kappa\x1)^*\cN$ is
trivial too.  So, since the second diagram above is commutative,
$(\vf\x1\x1)^*(A_C\x1)^*\cN$ is trivial.  Hence Equation~\Cs1.1) implies
that $(A_C\x1)^*\cN$ is trivial.

Fix isomorphisms,
        $$u\:(A_C\x1)^*\cN\risom\cO_{C\x J^1_C\x T}
        \and v\:(\kappa\x1)^*\cN\risom\cO_{P\x T},$$
 and set $\dr u:=(\Lambda\x1)^*v$.  Since the second diagram above is
commutative, $\dr u$ and
 $(\vf\x1\x1)^*u$
 differ by multiplication with an invertible (regular) function on $\dr
C\x J^1_C\x T$.  Since $\dr C$ is complete and integral, this function
is the pullback of an invertible (regular) function on $J^1_C\x T$.
Modifying $u$ accordingly, we may assume that $\dr u=(\vf\x1\x1)^*u$.

Set $R:=(C\x J^1_C)\x_{\CJ_C}(C\x J^1_C)$ and $\dr R:=R\x_{\CJ_C}P$, and
let $g\:\dr R\to R$ be the projection.  By
 Corollary~(5.5) in
 \cite{EGK},
 the second square in \Cs1.2) is Cartesian. Hence
$\dr R=(\dr C\x J^1_C)\x_P(\dr C\x J^1_C)$, and all
 three squares in the following diagram are Cartesian:
        $$\CD
        \dr R\x T @''' \dr C\x J^1_C\x T @>\Lambda\x1>>   P\x T \\
        @VVg\x1V @(\Fiberbox) @VV\vf\x1\x1V @(\fiberbox) @VV\kappa\x1V\\
          R\x T   @'''   C\x J^1_C\x T      @>A_C\x1>>   \CJ_C\x T
        \endCD$$

Form the two pullbacks $u_1,\ u_2$ of $u$ to $R\x T$ and
correspondingly those $\dr u_1,\ \dr u_2$ of $\dr u$ to $\dr R\x T$.
Now, $\dr u:= (\Lambda\x1)^*v$; so $\dr u_1=\dr u_2$.

The Abel map $A_C$ is smooth; see (3.1).  So the two projections from
$R\x T$ to $C\x J^1\x T$ are smooth too.  Hence the associated points of
$R\x T$ map to simple points of $C$.  Now, $\vf$ is an isomorphism off
the double point $Q$.
 Hence $g\:\dr R\to R$ is an isomorphism over the associated points of
$R$.  Since
        $$(g\x1)^*u_1=\dr u_1=\dr u_2=(g\x1)^*u_2,$$
 therefore $u_1=u_2$ holds at every associated point of $R\x T$, so
everywhere.

Since $u_1=u_2$, by descent theory $u$ descends to a trivialization of
$\cN$ on the image $V\x T$ of $A_C\x1$.  Now, $\CJ_C$ is a local
complete intersection of dimension $g$ by \cite{AIK76, (9), \p.8}, where
$g$ is the arithmetic genus of $C$.  Since each singular point of $C$ is
a double point,
 $\cod(\CJ_C-V,\, \CJ_C)\ge2$
 by Corollary~(6.8) in \cite{EGK}.
Hence, $\cN$ is the direct image of its restriction to $V\x T$.
Therefore, $\cN$ is trivial.  The proof is now complete in the case
where $\phi(T)\subset U^0$.  Call this the ``first case.''

Using our work in the first case, we will now establish the general
case.  To do so, fix an arbitrary rational point $t_0\in T$.  Let
$\cN_0$ be the fiber $\cN(t_0)$ viewed on $\CJ$.  We will prove that
$\cN_0$ is trivial.  Then we may conclude that $\phi(T)\subset U^0$, and
so we will have a complete proof of the autoduality theorem of (2.1).

By hypothesis, $\cN_0$ corresponds to a point of $U$.  So some multiple
$\cN_0^n$ corresponds to a point of $U^0$.  Moreover,
        $$A_\cL^*(\cN_0^n)=(A_\cL^*\cN_0)^n=\cO_C.$$
 Hence, by the preceding case with $S$ for $T$ and $\cN_0^n$ for $\cN$,
we may conclude that $\cN_0^n$ is the pullback of a sheaf on $S$.  Since
$S$ is a point, $\cN_0^n$ is trivial on $\CJ$.

Consider $A^*\cN_0$ on $C\x J^1$.  It defines a map
$\psi\:J^1 \to J$ such that $\psi[\cL]=0$.  Now, $\cN_0^n$ is trivial.
So $\psi(J^1)$ lies in the kernel of the $n$th power map $J\to J$.  This
kernel is finite.  Hence $\psi(J^1)=\{0\}$ since $\psi(J^1)$ is
connected and it contains 0.  Hence $\psi$ is constant since $J^1$ is
reduced.  Therefore, $A^*\cN_0$ is equal to the pullback of
some invertible sheaf $\cR$ on $J^1$.

Proceed by induction on $\delta$ as in the first case, but with $S$ for
$T$ and $\cN_0$ for $\cN$.  If $\delta =0$, then $U^0=U$ as we saw
above, and so $\cN_0$ is trivial by the first case.  If $\delta\ge1$,
then the argument in the first case goes through exactly as before since
the analogue of Equation~\Cs1.1) holds.  Thus $\cN_0$ is trivial, and
the proof is now complete.

 \dfn2 Let $C/S$ be a flat projective family of integral curves, $\cM$
an invertible sheaf of degree $m$ on $C/S$.  Define the {\it
translation\/} by $\cM$ to be the map,
        $$\tau_\cM\:\CJ^n\to\CJ^{m+n},$$
 given by tensoring  $\cM$ with a torsion-free sheaf.

More precisely, an $S$-map $t\:T\to\CJ^n$ corresponds to a torsion-free
rank-1 sheaf $\cN'$ on $C\x T'/T'$ with degree-$n$ fibers, where $T'/T$
is an \'etale covering, such that the two pullbacks of $\cN'$ to $C\x
T''$ are equal, where $T''\big/\, T'\x_TT'$ is an \'etale covering.  Let
$\cM'$ be the pullback of $\cM$ to $C\x T'$.  Then $\cM'\ox\cN'$ is a
torsion-free rank-1 sheaf of degree $m+n$ on $C\x T'/T'$, and its two
pullbacks to $C\x T''$ are equal.  Hence $\cM'\ox\cN'$ defines a map
$\tau_\cM(t)\:T\to\CJ^{m+n}$.

 \cor3 Let $C/S$ be a flat projective family of integral curves, $m$ and
$n$ integers, and $\cM$ an invertible sheaf of degree $m$ on $C/S$.  If
the curves only have double points
as singularities, then the
translation map $\tau_\cM$ induces an isomorphism,
        $$\tau_\cM^*\:\Pic^0_{\CJ^{m+n}/S}\risom\Pic^0_{\CJ^n/S},$$
 which is independent of the choice of $\cM$.  In particular, if $m=0$,
then $\tau_\cM^*$ is equal to the identity on $\Pic^0_{\CJ^n/S}$.
 \pf Note that $\tau_{\cO_C}=1_{\CJ^n}$.  And, if $\cM_1$ is also an
invertible sheaf on $C$, then
        $$\tau_\cM\circ\tau_{\cM_1}=\tau_{\cM\ox\cM_1}.$$
  So $\tau_\cM$ is an isomorphism, whose inverse is $\tau_{\cM^{-1}}$.
Hence $\tau_\cM^*$ is an isomorphism.  Moreover, if $\cM_1$ is of degree
$m$ too, then $\cM\ox\cM_1^{-1}$ is of degree 0, and it suffices to
prove that $\tau^*_{\cM\ox\cM_1^{-1}}=1.$   Thus we may assume $m=0$.

To prove that $\tau_\cM^*=1$, we may change the base by an \'etale
covering, and so assume that the smooth locus of $C/S$ admits a section
$\sigma$.  Set $\cL:=\cO_C(\sigma(S))$.  Then $\cL$  is an invertible
sheaf on $C$.  So,
  $$\tau_\cM=\tau_{\cL^{\ox n}}\circ\tau_\cM\circ\tau_{\cL^{\ox -n}}.$$
 Hence, since $\cL$ is of degree 1 on $C/S$, we may assume $n=0$.

Note that $\tau_\cM\circ A_\cL=A_{\cM\ox\cL}$.  Now,
$A_{\cM\ox\cL}^*=A_\cL^*$ by Proposition~(3.7).  Since
$A_\cL^*$ is an isomorphism by the autoduality theorem,
$\tau_\cM^*=1$, and the proof is complete.

 \cor4 Let $C/S$ be a flat projective family of integral curves.  If the
curves only have double points
as singularities, then the autoduality
isomorphism $\Pic^0_{\CJ/S}\risom J$ extends to a map
 $\eta\:\?U\to\CJ$,
 where $\?U$ is the natural compactification of $\Pic^0_{\CJ/S}$.
 \pf Since the curves have surficial singularities,
the projective
$S$-scheme $\CJ$ is flat, and its geometric fibers are integral; see
\cite{AIK76, (9), \p.8}.  Hence, by \cite{AK79II, Thm.~(3.1), \p.28},
there exists an $S$-scheme $U^=$ that parameterizes the torsion-free
rank-1 sheaves on the fibers of $\CJ/S$; the connected components of
$U^=$ are proper over $S$.  Moreover, $U^=$ contains $\Pic^0_{\CJ/S}$
as an open subscheme, and its scheme-theoretic closure in $U^=$ is, by
definition, $\?U$.  Furthermore, since $J/S$ is smooth and admits a
section (for example, the 0-section), by \cite{AK79II, Thm.~(3.4)(iii),
\p.40}, $\CJ\x\?U\big/\?U$ carries a universal sheaf $\cP$, which is
determined up to tensor product with the pullback of an invertible sheaf
on $\?U$.

The extension $\eta$ of the autoduality map is unique, if it exists.
Hence, by descent theory, it suffices to construct $\eta$ after
changing the base via an \'etale covering.  So we may assume that the
smooth locus of $C/S$ admits a section $\sigma$.  Set
$\cL:=\cO_C(\sigma(S))$.  Then $\cL$ is an invertible sheaf of degree 1
on $C/S$.  So the autoduality isomorphism is simply $A_\cL^*$, and
it suffices to prove that $(A_\cL\x1_{\?U})^*\cP$ is a torsion-free
rank-1 sheaf on $C\x\?U\big/\?U$.

The Abel map $A\:C\x J^1\to\CJ$ is smooth; see (3.1).  Hence
$(A\x1_{\?U})^*\cP$ is a torsion-free rank-1 sheaf on $C\x
J^1\x\?U\big/\?U$.  It suffices to prove that this sheaf is torsion-free
rank-1 on $C\x J^1\x\?U\big/(J^1\x\?U)$.  The sheaf is flat over
$J^1\x\?U$, by the local criterion, if its fiber is flat over the fiber
$J^1(u)$ for each $u\in\?U$.  Fix a $u$.  Making a suitable faithfully
flat base change $S'/S$, we may assume that the residue field of $u$ is
equal to that of its image in $S$.  Set $\cI:=\cP(u)$.  It suffices to
prove that $A(u)^*\cI$ is a torsion-free rank-1 sheaf on $C(u)\x
J^1(u)\big/J^1(u)$.

Suppose given an invertible sheaf $\cM$ of degree 0 on $C/S$.  Then the
translation map $\tau_\cM$ gives rise to the following commutative
diagram:
        $$\CD
        C\x J^1 @>A>> \CJ \\
        @V1_C\x\tau_\cM VV        @V\tau_\cM VV \\
        C\x J^1 @>A>> \CJ
        \endCD$$
 By Corollary~\Cs3), $\tau_\cM^*$ is the identity on $\Pic^0_{\CJ/S}$,
so on its closure $\?U$ too.  Thus
$\tau_\cM(u)^*\cI=\cI$.  Now, the diagram is commutative; hence,
        $$(1_C\x\tau_\cM(u))^*A(u)^*\cI=A(u)^*\cI.\eqno\Cs4.1)$$

Since $J^1(u)$ is integral, the lemma of general flatness applies, and
it implies that there is a dense open subset $W$ of $J^1(u)$ over which
$A(u)^*\cI$ is flat.  Now, by Part (ii)(a) of Lemma (5.12) on \p.85 of
\cite{AK80}, it is an open condition on the base for a flat family of
sheaves to be torsion-free rank-1 provided they are supported on a
family whose geometric fibers are integral of the same dimension.
Hence, since $A(u)^*\cI$ is torsion-free and of rank 1, after shrinking
$W$, we may assume that the restriction of $A(u)^*\cI$ to $C\x W/W$ is
torsion-free rank-1.  Fix a point $j_1$ of $W$ and an arbitrary point
$j_2$ of $J^1(u)$.

Making a suitable faithfully flat base change $S'/S$, we may assume that
each of $j_1$ and $j_2$ lies in the image of a section of $J^1/S$.
These sections represent invertible sheaves $\cM_1$ and $\cM_2$ of
degree 1 on $C/S$; set $\cM:=\cM_1\ox \cM_2^{-1}$.  Equation \Cs4.1)
implies that $A(u)^*\cI$ is torsion-free rank-1 over $\tau_\cM(u)^{-1}W$
as well.  Now, $j_2$ is an arbitrary point of $J^1(u)$.  Hence
$A(u)^*\cI$ is torsion-free rank-1 on $C(u)\x J^1(u)\big/J^1(u)$, and
the proof is complete.

 \references
 \serial{adv}{Adv. Math.}
 \serial{ajm}{Amer. J. Math.}
 \serial{bams}{Bull. Amer. Math. Soc.}
 \serial{bpm}{Birkh\"auser Prog. Math.}
 \serial{ca}{Comm. Alg.}
 \serial{crasp}{C. R. Acad. Sci. Paris}
 \serial{MathScand}{Math. Scand.}
 \def\GrothFest#1 #2 {\unskip, \rm
 in ``The Grothendieck Festschrift. Vol.~#1,'' P.~Cartier, L.~Illusie,
N.M.~Katz, G.~Laumon, Y.~Manin, K.A.~Ribet (eds.) \bpm 86 1990 #2 }
 \def\osloii{\unskip, \rm in ``Real and complex singularities. Proceedings,
Oslo 1976,'' P. Holm (ed.), Sijthoff \& Noordhoff, 1977, pp.~}
 \def\sitgesii{\unskip, \rm in ``Enumerative Geometry. Proceedings,
Sitges, 1987,'' S. Xamb\'o Descamps (ed.), \splnm 1436 1990 }
 \def\splnm #1 #2 {Springer LNM {\bf#1}\ (#2)}

AIK76
 {A. Altman, A. Iarrobino, and S. Kleiman},
 Irreducibility of the compactified Jacobian
 \osloii 1--12.

AK76
 A. Altman and S. Kleiman,
 Compactifying the Jacobian,
 \bams 82(6) 1976 947--49

AK79
 A. Altman and S. Kleiman,
 Bertini theorems for hypersurface sections contaning a subscheme,
 \ca 7(8) 1979 775--790

AK79II
 A. Altman and S. Kleiman,
 Compactifying the Picard scheme II,
 \ajm 101 1979 10--41

AK80
 A. Altman and S. Kleiman,
 Compactifying the Picard scheme,
 \adv 35 1980 50--112

AK90
 A. Altman and S. Kleiman,
 The presentation functor and the compactified Jacobian
 \GrothFest I 15--32

E95
 D. Eisenbud,
 Commutative Algebra,
 Springer GTM {\bf 150} (1995).

E99
 E. Esteves,
 Compactifying the relative Jacobian over families of reduced curves,
 to appear in Trans. Amer. Math. Soc..

EGK
 E. Esteves{,} M. Gagn\'e{,} and S. Kleiman,
 Abel maps and presentation schemes,
 to appear.

FGA
 A. Grothendieck,
 Technique de descente et th\'eor\`emes d'existence en g\'eom\'etrie
alg\'e\-bri\-que V, S\'eminaire Bourbaki, 232, Feb. 1962, and VI,
S\'eminaire Bourbaki, 236, May 1962.

EGA
 A. Grothendieck and J. Dieudonn\'e,
 El\'ements de G\'eom\'etrie Alg\'ebrique,
 Publ. Math. IHES, O$_{\rm III}$, Vol. 11, 1961;
III$_2$, Vol. 17, 1963; IV$_3$, Vol. 28, 1966; and
IV$_4$, Vol. 32, 1967.

K87
 S. Kleiman,
 Multiple-point formulas II: the Hilbert scheme
 \sitgesii , 101--138.

KM76
 F. Knudsen and D. Mumford,
 The projectivity of the moduli space of curves I,
 \MathScand 39 1976 19--55

L59
 S. Lang,
 Abelian varieties,
 Interscience Tract {\bf 7} (1959).

M60
 A. Micali,
 Alg\`ebre sym\'etrique d'un id\'eal,
 \crasp 251 1960 1954--56

M65
 D. Mumford,
 Geometric invariant theory,
 Springer Ergebnisse {\bf 34} (1965).

M70
 D. Mumford,
 Abelian Varieties,
 Oxford Univ. Press, 1970.

R67
 M. Raynaud,
 Un th\'eor\`eme de representabilit\'e relative sur le foncteur de Picard,
 Expos\'e XII, in ``Th\'eorie des intersections et th\'eor\`eme de
Riemann--Roch (SGA 6),'' Springer LNM {\bf 225} (1971).

\endreferences
\mark{E}\bye